\documentclass{article}

\usepackage[preprint]{corl_2024} 

\usepackage{amsmath}

\DeclareMathOperator*{\argmin}{arg\,min}

\usepackage{amsfonts}
\usepackage{graphicx}
\usepackage{amssymb}
\usepackage{caption} 
\usepackage{tabularx}
\usepackage{subcaption}
\usepackage{bm}
\usepackage{algorithm2e}
\usepackage{booktabs}
\usepackage{geometry} 
\usepackage{siunitx} 
\newcommand{\norm}[1]{\left\lVert#1\right\rVert}
\usepackage{multirow}
\newcommand{\tr}{\mathsf{T}}  
\newcommand{\Rb}{\mathbb{R}} 
\usepackage{tikz}
\usepackage{soul}
\usepackage{bbding}
\usepackage[export]{adjustbox}
\usepackage{hhline}
\usepackage{makecell}
\usepackage{caption}
\usepackage{blindtext}
\usepackage{float}
\usepackage[export]{adjustbox}
\usepackage{bbm}

\title{Second-Order Stein Variational Dynamic Optimization}
 
\author{
  Yuichiro Aoyama\\
  \texttt{yaoyama3@gatech.edu}
  \And
  Peter Lehmann \\
  \texttt{plehmann7@gatech.edu}
  \And
  Evangelos A. Theodorou \\
  \texttt{evangelos.theodorou@gatech.edu}\\
  \\
  School of Aerospace Engineering\\
  Georgia Institute of Technology 
  United States\\
}

\begin{document}
\maketitle


\begin{abstract}
    We present a novel second-order trajectory optimization algorithm based on Stein Variational Newton's Method and Maximum Entropy Differential Dynamic Programming. The proposed algorithm, called Stein Variational Differential Dynamic Programming, is a kernel-based extension of Maximum Entropy Differential Dynamic Programming that combines the best of the two worlds of sampling-based and gradient-based optimization. The resulting algorithm avoids known drawbacks of gradient-based dynamic optimization in terms of getting stuck to local minima, while it overcomes limitations of sampling-based stochastic optimization in terms of introducing undesirable stochasticity when applied in online fashion. To test the efficacy of the proposed algorithm, experiments are performed for both trajectory optimization and model predictive control. The experiments include comparisons with unimodal and multimodal Maximum Entropy Differential Dynamic Programming as well as Model Predictive Path Integral Control and its multimodal and Stein Variational extensions. The results demonstrate the superior performance of the proposed algorithms and confirm the hypothesis that there is a middle ground between sampling and gradient-based optimization that is indeed beneficial for the purposes of dynamic optimization. This middle ground consists of different mechanisms that combine sampling with gradient-based optimization. In this paper, we investigate these different mechanisms and show their benefits in dealing with non-convex dynamic optimization problems found in trajectory optimization and model predictive control.
    
    

\end{abstract}

\keywords{Optimization, Variational Inference, Robots} 


\section{Introduction}

   Trajectory Optimization (TO) is essential for robotic and autonomous systems to operate safely in dense and cluttered environments. From an optimization point of view, the planning and trajectory optimization problems in robotics have some unique characteristics. These include the non-convexity of the objective function in consideration, the nonlinear state and actuation constraints, and the highly nonlinear dynamics that characterize most robotic systems. The characteristic of non-convexity is due to the existence of multiple local minima. Such situations arise in reaching and navigation tasks performed in cluttered and dense environments.

  There is a long history of TO methods and their application to robotics and autonomous systems. Starting from the  Differential Dynamic Programming (DDP) \citep{Jacobson1970ddp} and the Iterative Pontryagin Maximum Principle (iPMP) \citep{Pontryagin1964mathematical}, these optimal control methods have been extensively used in robotics. Both DDP and iPMP belong to the same family of algorithms that are gradient-based and rely on forward-backward sweeps. Their mathematical structure is derived using either the Bellman Principle of Optimality or the Pontryagin Maximum Principle and consists of the backpropagation of the value function or costate and forward propagation of state trajectories. This step of backpropagation is performed locally around the trajectories that are updated at each iteration. While the locality of the backward pass makes DDP and iPMP scalable and fast, it often results in sub-optimal solutions. 
  
  

  On the opposite side of gradient-based trajectory optimization, sampling-based stochastic optimization methods overcome the issue of locality and have no differentiability requirements. These algorithms have been successfully used in robotics for different systems and tasks and are typically deployed in Model Predictive Control (MPC) fashion.  The most representative algorithms in this area are the Cross-Entropy (CE) as well as Model Predictive Path Integral control (MPPI) and its variations \cite{Williams2016MPPI, Williams2018MPPI,RobustMPPI,Williams-RSS-18,Wang2021Tsallis}. In general, these algorithms have fewer requirements and assumptions compared to more traditional approaches to dynamic optimization since dynamics and cost functions do not have to be differentiable. This characteristic makes these methods easy to use. In addition, the sampling process in stochastic optimization provides exploration capabilities and results in optimization steps that are less incremental and more exploratory than gradient-based optimization. The drawback of MPPI is that the resulting optimal decisions are noisy, and this can result in undesirable stochastic behaviors. This characteristic is attributed to the added stochasticity required for sampling trajectories. Recent efforts to reduce the effects of added stochasticity include sampling from colored noise \citep{Vlahov2024ColoredMPPI} and sampling using the inverse of the Hessian of the cost function in consideration \citep{yi2024covompc}. While both works provide meaningful ways to address the issue, they impose extra requirements related to the forms of sampling and the need to have a twice differentiable cost function. 
  
While the sampling aspect of stochastic optimization is a way to handle complex objective functions, this feature does not completely address the issue. Naive sampling quite often can be inefficient in exploring the state space. These limitations motivated the work on multi-modality. In  \cite{urain2022learaningimplicitprior}, an energy based model that represents multimodality is learned and incorporated in gradient- and stochastic-based optimization for planning. Reference
\cite{Carvalho2023Diffusionrobot} proposes a diffusion model to handle large datasets with multimodality to learn prior of a trajectory optimization task, while the work in 
 \cite{le2023acceleratingOT} introduces an Optimal Transport technique to plan trajectories for a task with a multimodal objective. 
 Although the aforementioned approaches advance planning, they have a significant limitation in that they rely in a linear time invariant dynamics representations. Alternative algorithms \cite{lambdert2021SVMPC,lambert2021entropyregularizedmotionplanning,power2024constSV} explore multi-modality by making use of Stein Variational Gradient Descent (SVGD) \cite{Liu2017SteinGD}. The work in  \cite{lambdert2021SVMPC} is a Stein Varational extensions of MPPI while the work in \cite{lambert2021entropyregularizedmotionplanning} is only used for planning tasks. Finally reference \cite{power2024constSV} proposes a Stein Variational constrained trajectory optimization algorithm. The approach scales unfavorably with the dimensionality of control, state variables, time horizon and number of equality/inequality constraints including dynamics. To alleviate this issue short time horizons are chosen when deployed in MPC fashion resulting in myopic behaviors.    
  

   In this paper, we investigate two mechanisms that combine sampling with gradient-based optimization and parallelism, and develop algorithms for planning and MPC with highly non-convex objectives.  The first mechanism is a schedule of optimization steps, that consists of number of Netwon steps followed by an iteration that relies on sampling.  The motivation for this hybrid mechanism is that while Newton  steps  can locally optimize the objective function,  sampling can disrupt this process to ensure that optimization is not getting trap into a sub-optimal local minima.  In the context of TO for  robotics, this mechanism was first introduced in the work on Maximum Entropy DDP (MEDDP) and its unimodal and multimodal extensions. The difference of this work with respect to \cite{so2022maximum} is that here we show the benefits of this mechanisms in the context of MPC. 
    It is also interesting to note that this hybrid optimization mechanism, which we name \textit{ Messy Optimization}, was also proposed in the area of business administration and management science \cite{harford2016messy} as mental model for enhances creativity and improving innovation. 

    The second mechanism goes one step further than MEDDP and results in new algorithm, the Stein Variational Differential Dynamic Programming (SVDDP). Grounded on the Stein Variational Newton's Method (SVNM) \cite{Chen2019ProjectedSVNewton}, SVDDP is a kernel-based extension of MEDDP. In SVDDP, the optimization process alternates between two phases.
    Starting from different initial trajectories/particles, the first phase consists of multiple DDP updates  performed in parallel fashion resulting in optimized state and control trajectories. In the second phase, these trajectories are  updated only once with DDP using the SVNM to update the controls and generate the new state trajectories. This phase enforces exploration by pushing the trajectories to stay away from each other. In addition, the fact that this update is performed once every few iterations of the first phase ensures that the overall optimization will not be disrupted.   We benchmark the proposed algorithm against variations on MEDDP and MPPI for both TO and MPC and show that it matches and, in most cases, outperforms these algorithms. The rest of the paper is organized as follows. Section \ref{sec:background} introduces the background of our work which includes MEDDP and SVGD and SVNM. In Section \ref{sec:SVDDP} we derive SVDDP. Section \ref{sec:result} includes experimental results, and section \ref{sec:conclusion} concludes with future extensions.

 

\section{Background}\label{sec:background}
\subsection{Maximum Entropy Dynamic Second-Order c Optimization}
Consider a TO problem for a dynamical system with state $x\in \mathbb{R}^{n_{x}}$ and control $u\in \mathbb{R}^{n_{u}}$. Given the deterministic dynamics
$x_{t+1} = f(x_{t}, u_{t})$,
The TO problem is formulated as follows:   
\begin{equation}\label{eq:ocp_cost}
\min_{X,U} J(X,U) = \min_{X,U} \Phi(x_{T}) +  \sum_{t=0}^{T-1}l_{t}(x_{t},u_{t}), \quad \text{s.t.}  \quad x_{t+1} = f(x_{t}, u_{t}).
\end{equation}
There exist several second-order solvers that can solve the problem above. They solve Quadratic Programming (QP) with a quadratic approximation of cost under constraints arising from linearized dynamics \cite{Gill2002SNOPT, Wachter2006IPOPT}. However, relying on local information on the cost and dynamics, they are vulnerable to being captured at poor local minima. 
To promote exploration while solving the problem, we consider a stochastic control policy $\Pi(U) = \Pi(u_{0}, \cdots, u_{T-1})$, and 
introduce an entropy term $H[\Pi] = - \mathbb{E}_{\Pi}[\ln (\Pi(U))]$ to the objective. Considering the expectation, we have a new objective:
\begin{equation}\label{eq:ME_DDP_objective}
J_{\Pi}(x_{0}, \Pi)= \mathbb{E}_{\Pi} \big[ \sum_{t=0}^{T-1}l(x_{t},u_{t})  + \Phi(x_{T})\big] - \alpha H[\Pi],
\end{equation}
with $\int \Pi(U) {\mathrm{d}}U= 1$. The term $\alpha$ is the temperature that governs the relative importance of the entropy in the new objective.  With this new objective, the  optimal control problem seeks a stochastic policy that minimizes the initial deterministic cost function while maximizing the corresponding Shannon entropy. The solution of the minimization of the new objective is given by 
\begin{align}
\Pi^{\ast}(U) =Z_{0}^{-1} \exp[-J(x_{0}, U)/\alpha],   
\end{align}
where $Z_{0}$ is a partition function. We refer the reader to Appendix \ref{sec:app_DDP} for the derivation. 
One way to efficiently solve this problem is to use Differential Dynamic Programming (DDP).



\subsection{Maximum Entropy DDP}
\textbf{Unimodal Case:} Application of the Bellman's principle in minimizing \eqref{eq:ME_DDP_objective} results in the optimal policy $\pi_{t}^{\ast}$:
\begin{equation} \label{eq:me_ddp_pi_star_uni}
\pi^{\ast}(u|x) = {Z(x)}^{-1}\exp[ -Q(x,u)/{\alpha}], \quad \text{with}\ 
V(x) = -\alpha \ln Z(x),
\end{equation}
where $ Q(x,u) $ is the state-action value function ($Q$ function) of DDP and $Z(x)$ is the corresponding partition function. We refer readers to \citep{so2022maximum} and Appendix \ref{sec:app_DDP} for derivation. This indicates that the optimal policy has a form of unimodal Gaussian $\pi^{\ast}(\delta u|\delta x) \sim \mathcal{N}(\delta u^{\ast}, \alpha Q_{uu}^{-1})$, where $\delta u^{\ast}$ is a solution of vanilla deterministic DDP:
\begin{equation}\label{eq:delta-u-star-main}
\delta {u}^{\ast}= \kappa +K\delta x, \quad
\text{with}\ \kappa = -{Q}^{-1}_{{uu}}{Q_{{u}}},\ K = -{Q}^{-1}_{{uu}}{Q_{{ux}}}.
\end{equation}
\textbf{Multimodal Case:} In the multimodal case, we consider the LogSumExp approximation of the objective given $N$ trajectories, or modes. The exponential transformation $
\mathcal{E}_{\alpha}(y) = \exp[
-y/{\alpha}]$ of the value function results in a control policy represented as a mixture of Gaussians whose categorical distribution is proportional to the value function of each trajectory. This multimodal policy is represented as follows: 
\begin{align}\label{eq:me_ddp_pi_star_multi}
 \pi^{\ast}(u|x) 
= \sum_{n=1}^{N} \omega^{(n)}(x) \pi^{\ast (n)}(u|x) \quad 
\text{with} \quad \omega^{(n)} 
\propto 
\mathcal{E}_{\alpha}[
 V^{(n)}(x)], \quad \text{and} \quad \sum_{n=1}^{N}\omega^{(n)} =1. 
\end{align}
We note that these modes may collapse to a small number during optimization. This phenomenon reduces the capability of multimodal MEDDP to explore the state space and avoid local optima.

\subsection{Stein Variational Gradient Descend and Newton's method}
This section briefly covers Stein Variational Gradient Descent (SVGD) and SV Newton's Method (SVNM). We first introduce the update law of SVGD, then describe the connection to functional gradient descent. Finally, we review SVNM as an extension to the second-order optimization of SVGD. For further information, we refer readers to \citep{Liu2017SteinGD} for SVGD and \citep{detommaso2018steinNewton} for SVNM. In this section, we use $x \in \mathbb{R}^{d}$ to denote a variable of a function but not a state of a system. In addition, we use $x_{n} \in \mathbb{R}^{d}$ to denote a sampled point.

\textbf{Stein Variational Gradient Descend:} 
Let $p$ on $\Rb^{d}$ be a target distribution that we would like to approximate using a collection of samples. We drop the argument of the function $p$, and specify it only when we evaluate it at point $x_{n}$ as $p(x_{n})$. With samples $\{x_{n}\}$ from a tractable reference distribution $q$ on $\Rb^{d}$, SVGD iteratively computes a transport map $T:\mathbb{R}^{d} \rightarrow \mathbb{R}^{d}$ so that the transformed samples from $q$, i.e., $\{T(x_{n})\}$ can be an empirical approximation of $p$. This map is obtained by solving the following optimization problem.
\begin{align}\label{eq:SVGD_problem}
\min_{T} \{ D_{\rm{KL}}( q^{l+1}||p) = \mathbb{E}_{T\ast q^{l}}\big[\log( q^{l+1}/p)\big] \}, \quad \text{with} \ T\ast q^{l} = q^{l+1},
\end{align}
where $l$ stands for $l$ th iteration. Here, $D_{\rm{KL}}$ is KL divergence that measures the difference between the two distributions. To simplify this problem, SVGD considers the vector-valued Reproducing Hilbert Kernel Space (RKHS) $\mathcal{H}^{d} = \mathcal{H}\times \cdots \times \mathcal{H}$, where $\mathcal{H}$ is a scalar-valued RKHS with kernel $k(x, x')$. With RKHS, one can represent a function as a weighted composition of kernels, centered at samples. This representation also allows us to perform optimization in functional space by working with the weights in the same way as the normal optimization technique, such as GD. Here, the transportation map is written as a perturbation of the identity map denoted by $I$ along $\mathcal{H}^{d}$ as
\begin{align*}
T^{l}(x) = I(x) + Q(x), \quad \text{for} \quad Q(x) \in \mathcal{H}^{d}.
\end{align*}
 With this $T^{l}$, we define the objective of the problem \eqref{eq:SVGD_problem} as
\begin{align}
\hat{J}_{q^{l}}[Q] = D_{\rm{KL}}\big((I + Q)\ast q^{l} || p\big ),
\end{align}
which should be decreased by the proper choice of $Q$ from $\hat{J}_{q^{l}}[\bm{0}]$, where $\bm{0}(x)=0$ is a zero map. Since this $Q$ lies in the steepest direction of the functional $\hat{J}_{q}^{l}$ at $\bm{0}$, the map is given as follows. 
\begin{align}\label{eq:SVGD_T}
T^{l} = I - \epsilon \nabla \hat{J}_{q^{l}}[\bm{0}],
\end{align}
with a step size $\epsilon$. 
Observe the similarity with GD update of a normal optimization, i.e., $x^{l+1} = x^{l} + \epsilon \nabla f_{0}(x^{l})$ with objective $f_{0}$.
The descent direction satisfies the following condition with $S\in\mathcal{H}^{d}$,
\begin{align*}
D\hat{J}_{q}[S](V) = \langle \nabla \hat{J}_{q(x)}[S], V \rangle_{\mathcal{H}^{d}}, \quad \forall V \in \mathcal{H}^{d}, 
\end{align*}
where we drop $l$. Here, the left-hand side is the first variation of $\hat{J}_{q}$ at $S$ along $V$ defined as:
\begin{equation*}
D\hat{J}_{q}[S](V) = \lim_{\tau \rightarrow 0} \frac{1}{\tau}(\hat{J}_{q}[S+\tau V] - \hat{J}_{q}[S]).
\end{equation*}
\citet{Liu2017SteinGD} showed that the gradient at zero has a closed-form solution and is empirically approximated by $N$ particles as:  
\begin{equation}\label{eq:SVGD_emp_func_gradient}
-\nabla \hat{J}_{q}[\bm{0}](x_{i}) = \frac{1}{N}\sum_{n=1}^{N}
\Big[[\nabla_{x_{n}}\log p(x_{n})]k(x_{n}, x_{i}) + \nabla_{x_{n}}k(x_{n}, x_{i})\Big] = \phi^{\ast}(x_{i}),
\end{equation}
where the first term in the summation follows the gradient direction, and the second term spreads the particles apart from each other and is thus known as the repulsive force term. For a particle, the update law is given with a step size $\epsilon$ as follows:
\begin{equation}\label{eq:SV_update}
x_{n}^{l+1} \gets x_{n}^{l} + \epsilon \phi^{\ast}(x_{n}^{l}).
\end{equation}
\textbf{Stein Variational Newton's Method:} To perform the second-order optimization in functional space, we define the second variation of $\hat{J}_{q}$ at $\bm{0}$ along the pair of directions $V, W \in \mathcal{H}^{d}$ as:  
\begin{align*}
D^{2}\hat{J}_{q}[\bm{0}](V,W) = \lim_{\tau \rightarrow 0} \frac{1}{\tau}(D \hat{J}_{q}[\tau W](V) - D\hat{J}_{q}[\bm{0}](V)).
\end{align*}
The Newton direction $W$ is obtained by the optimality condition given by the following equation. 
\begin{align*}
D^{2}\hat{J}_{q}[\bm{0}](V,W) =  - D \hat{J}_{q}[\bm{0}](V), \quad \forall V \in \mathcal{H}^{d}, 
\end{align*}
which gives the transformation map of the Newton direction as a perturbation of the map $I$ as: 
\begin{align*}
T = I + \epsilon W.
\end{align*}
\citet{detommaso2018steinNewton} proved that the Newton direction $W = (w_{1}, \cdots, w_{d})^{\mathsf{T}}$ satisfies 
\begin{align}\label{eq:SV_Newton_hessian}
&\sum_{i=1}^{d}\Bigl\langle\sum_{j=1}^{d} \langle h_{ij}(y,z),w_{j}(z) \rangle_{\mathcal{H}} + \nabla \hat{J}[\bm{0}]_{i}(y), v_{i}(y) \Bigr\rangle_{\mathcal{H}} = 0, \quad \forall \ V = (v_{1}, \cdots v_{d})^{\tr} \in \mathcal{H}^{d}, \\\notag
&\text{where} \quad 
h_{ij}(y,z) = \mathbb{E}_{x \sim q}\big[-[\nabla^{2}\log p(x)]_{i,j}k(x, y)k(x,z) + \nabla k{(x,y)}_{i}\nabla k(x,z)_{j}\big].
\end{align}
In the same work, the Galarkin approximation of the solution of $W$ is also proposed, where $W$ is expanded on $\mathcal{H}^{d} = \mathrm{span}\{k(x_{1}, \cdot), \cdots, k(x_{N}, \cdot)\}$. This expansion leads to the following approximation with unknown coefficients $\beta$ as 
\begin{align*}
[w(z)]_{j} = \sum_{n=1}^{N}[\beta^{n}]_{j}k(x_{n},z), \quad \text{for} \ j = 1, \cdots d.
\end{align*}
Here, the superscript of $\beta$ is an index that indicates samples. In addition, $j$ stands for $j$ th element of a vector. The $\beta$s are given as a solution of linear systems given by 
\begin{equation*}
\sum_{n=1}^{N} H^{s,n} \beta ^{n} = -\nabla \hat{J}_{p}[\bm{0}](x_{s}), \quad \text{for} \ s = 1, \cdots ,N,
\end{equation*}
with $\beta^{n} = (\beta^{n}_{1}, \cdots ,\beta^{n}_{d})$, and Hessian $H^{s,n}_{i,j} = h_{ij}(x_{s}, x_{k})$. The authors also propose the block diagonal approximation of the system to enjoy parallelization, which transforms the system as: 
\begin{align}\label{eq:SV_Newton_blkdiag_approx}
 H^{n,n} \beta ^{n} = -\nabla \hat{J}_{p}[\bm{0}](x_{n}), \quad \text{for} \ n = 1, \cdots, N ,
\end{align}
SV Newton's Method repeats the process of solving the linear systems above and updating particles
\begin{align*}
x_{n}^{l+1} \gets x_{n}^{l} + \epsilon w(x_{n}). 
\end{align*}

\section{Stein Variational Dynamic Second-Order Optimization}\label{sec:SVDDP}
In this section, we apply SVNM to DDP to derive a new algorithm that can explore multiple local minima as MEDDP. Specifically, we wish to use the repulsive force term to spread trajectories running in parallel so that the policies of multiple DDPs do not collapse. In the MEDDP setting, the optimal distribution is $\pi^{\ast}$ in \eqref{eq:me_ddp_pi_star_uni}. Therefore, we substitute $\pi^{\ast}(u)$ for $p(x)$ of the target distribution of SV methods, which leads to a repulsive force term for MEDDP. After the derivation, we explain three essential components of the algorithm: the kernel, the step size, and the heuristic about the frequency of applying the SV update. Finally, we explain how constraints are handled within DDP. We provide the pseudocode in Appendix \ref{sec:app_svddp_code}.

\subsection{ME DDP with Functional Gradient and Hessian}
\begin{figure}[!h]
\centering
\includegraphics[trim={2cm 16cm 4.8cm 2.5cm},clip,width=1.0\textwidth]{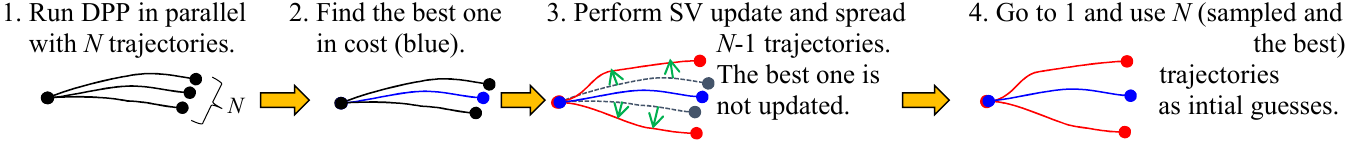}
\caption{Schematic of the SVDDP.}
\label{fig:SVDDP_alg}
\end{figure}

We consider $N$ trajectories optimized by DDP that we would like to spread out for exploration. Using the deviation $u^{\ast} = \bar{u} + \delta u^{\ast}$, the gradient operator $\nabla_{u\ast}$ for $Q$ function of DDP is equivalent to $\nabla_{\delta u}$. Therefore, with the quadratically approximated $Q$, and the substitution of $\pi^{\ast}(u^{\ast})$ for $p(x)$ in the functional gradient \eqref{eq:SVGD_emp_func_gradient}, we get: 
\begin{align}\label{eq:SVDDP_grad}
-\nabla \hat{J}_{q}[\bm{0}](u^{\ast(s)})&= 
\frac{1}{N}\sum_{n=1}^{N}
\Big[-\frac{1}{\alpha}\underbrace{(Q_{u}^{(n)}+ Q_{ux}^{(n)}\delta x^{\ast(n)} +  Q_{uu}^{(n)}\delta u^{\ast(n)})}_{=0} k(u^{\ast (s)}, u^{\ast (n)}) \\\notag
& \hspace{50mm} + \nabla_{u^{\ast( n)}}k(u^{\ast (n)}, u^{\ast (s)})\Big], \ \text{for} \ s = 1, \cdots,N,
\end{align}
where the derivatives of $Q$ are evaluated at $\bar{u}$.
The first term on the right-hand side is zero because of the optimality condition used in the derivation of DDP. Therefore, we are left with the repulsive force terms, which seem attractive for spreading trajectories. However, since we lose the effect of $\alpha$, this formulation cannot capture the relative significance of the entropy in the objective. We need SVNM to capture it. By substituting $\pi^{\ast}(u)$ for $p(x)$ in the approximated Hessian \eqref{eq:SV_Newton_blkdiag_approx} with the terms in \eqref{eq:SV_Newton_hessian}, and performing empirical approximation, the Hessian for SV Newton's method is obtained as: 
\begin{align}\label{eq:Hss_SVDDP}
H^{s,s} = \frac{1}{N}\sum_{n=1}^{N}\Big[ \frac{Q^{(n)}_{uu}}{\alpha} k (u^{ (n)},  u^{(s)})^{2} + \nabla_{ u^{(n)}} k(u^{(n)}, u^{(s)})\nabla_{u^{(n)}} k( u^{(n)},  u^{(s)})^{\tr}\Big],
\end{align}
where we drop $\ast$. By plugging this into \eqref{eq:SV_Newton_blkdiag_approx} and solving the equations, we get SV update of control for DDP
\begin{align}\label{eq:SV_ddp_update}
u^{(s)} \gets u^{(s)} +\epsilon w(u^{(s)}) + K \delta x, \quad \text{with} \quad w(u^{(s)}) = \sum_{n=1}^{N} \beta^{(n)}k( u^{(n)},  u^{(s)}),
\end{align}
where $K$ is feedback gain of DDP in \eqref{eq:delta-u-star}. 
Now, $\alpha$ is captured in $H^{s,s}$. Hence, we believe that this should be the correct formulation to incorporate the SV method with DDP. The second term in \eqref{eq:Hss_SVDDP} is known as an expected outer product of the gradient, which captures information on the variation in all directions \citep{Trivedi2014expectoutergradient}. Thus, it can be seen as an approximated Hessian of the kernel. We also provide another example of interpretation in the Appendix \ref{sec:app_outerproduct_kernel}. 
To comprehend how the Hessian affects the repulsive force term, let us consider an example with $H \in \Rb^{2 \times 2}$. We can draw an ellipsoid from its eigenvectors and eigenvalues whose major and minor axes point to the directions where $Q$ and $k$ increase most rapidly and slowly, respectively.  When solving \eqref{eq:SV_Newton_blkdiag_approx}, and $H$ is inverted, the scaling factors of the axes, i.e., eigenvalues, are also inverted. However, the eigenvectors are preserved. Consequently, the axes are flipped. Hence, the major and minor axes of $H^{-1}$ indicate the direction in which $Q$ and $k$ increase most slowly and rapidly, respectively. When \eqref{eq:SV_Newton_blkdiag_approx} is solved, the repulsive force term is projected in the direction where $Q$ and $k$ increases slowly. This means that the force is transformed so that it does not disturb the optimization process too 
much and spreads trajectories more by the curvature information on $Q$ and $k$, respectively. This projection leads to an efficient exploration.  We provide the schematic of SVDDP in Fig.\ref{fig:SVDDP_alg}. We note that the update \eqref{eq:SV_ddp_update} is not applied to the best trajectory, ensuring that at least one mode shows a monotonic cost reduction to preserve the convergence property of DDP as in \citep{so2022maximum, Dong2021Replicaexchange}. We note that we use cost as a metric, although other metrics e.g., diversity, can be used depending on the purpose.

\subsection{Essential Components in SVDDP}\label{sec:components_SVDDP}
\textbf{Kernel and Property of Repulsive Force:} The source of repulsive force in the SV method is the gradient of the kernel; thus, the choice of kernel is critical. As in other SV literature, we use the Gaussian kernel and choose its length parameter by median heuristic \citep{garreau2018large_median_heuristic}. Details of how it helps to generate repulsive force are provided in Appendix \ref{sec:app_kernel_repulsive_force}.

\textbf{Step Size:} Since the step size $\epsilon$ determines how much extent the trajectories are disturbed, its choice is important. In SV literature, \citep{Chen2019ProjectedSVNewton, Chen2020ProjectedSVGD} use backtracking line search. In our setting, however, because $w(u)$ only has information of repulsive force term, we cannot apply the same strategy. For robotic applications with SV, \citep{lambdert2021SVMPC, Wang2021Tsallis} use different fixed $\epsilon$s depending on tasks. During implementation, we observe that too large $\epsilon$ yields \texttt{NaN}s in the state of the dynamical system especially when it is unstable. Our strategy is to try with a fixed $\epsilon$ and if we encounter \texttt{NaN}, we back-track $\epsilon$ until we obtain a real number. Although this works well, the proper choice of $\epsilon$ remains future work.  

\textbf{Heuristics for Update Frequency:} In SVDDP, we allow DDP to perform deterministic iterations and then disrupt the process by conducting an update using SV. This is the same heuristic that has been used in MEDDP and its unimodal and multimodal variations \citep{so2022maximum}.  

\subsection{Handling Constraints}\label{sec:handle_const_DDP}
This section provides how we handle state and control constraints in DDPs. MEDDP was proposed in an unconstrained setting. It encouraged systems to keep distance from obstacles by adding penalty terms to the cost function \citep{so2022maximum}.
Our work use relaxed $\log$ barrier function to handle constraints, which is similar to \cite{Grandia2019BarrierDDP}, but with better conditioning of the Hessian with the Gauss-Newton approximation \cite{Nocedal2006numerical}. The exact barrier can always make the trajectories feasible, but cannot accept infeasible initial guesses \cite{Frisch1955logarithmic, Fiacco1968sequential}. This is problematic because the update in \eqref{eq:SV_update} and the sampling from \eqref{eq:me_ddp_pi_star_uni} and \eqref{eq:me_ddp_pi_star_multi} are not explicitly aware of the constraints, making the trajectories for the next DDP infeasible. With the relaxed barrier, DDPs can perform optimization even when sampled trajectories are infeasible. Although it is relaxed and thus the constraints can be violated, the barrier function can help quickly recover feasibility by imposing a high cost on infeasibility. Finally, we note that for exploration, the control limits are not considered to promote exploration. They are considered only in the optimization. The details of constrained DDP are in Appendix \ref{section:app_DDP_with_barrier}.

\section{Experimental Results}\label{sec:result}
This section provides experimental results, including TO and MPC. For TO, we compare normal DDP, Unimodal, Multimodal Gaussian MEDDPs (UG-, MG-MEDDP), and SVDDP. We note that DDPs are modified so that they can handle constraints as in section \ref{sec:handle_const_DDP}. For MPC, our comparisons involve DDP-MPC, UG-MEDDP, MG-MEDDP, SVDDP in MPC mode. In addition, MPPIs with three different policy parameterization and update laws are included in the comparison. They are Unimodal Gaussian MPPI (UG-MPPI), Multimodal Gaussian MPPI (MG-MPPI), and Unimodal Gaussian with Stein Variational update law (SV-MPPI) \citep{Wang2021Tsallis}. These three MPPIs correspond to the UG-MEDDP, MG-MEDDP, and SVDDP, respectively. In terms of systems, we use three systems, including a 2D Car, quadrotor, and 7DoF arm in TO and 2D car and quadrotor for MPC. In addition, we have comparisons of DDPs in swimmer and ant dynamics in Brax \cite{brax2021github} whose results are presented in section \ref{sec:add_ant_swimmer}.


\noindent \textbf{Trajectory Optimization:} This experiment is to examine the exploration capability and convergence of SV and ME DDPs over normal DDP. In addition, we observe the effect of temperature $\alpha$, number of modes, and stepsize $\epsilon$ on the performance. Figure \ref{fig:TO_overview} provides an overview of the experiments. We performed the same experiment 15 times to capture the properties of the algorithms, and 15 outcomes are shown in the figure. The task is to reach the target while avoiding obstacles under control limits. The configuration and shapes of the obstacles give multiple local minima, which normal DDP without exploration gets captured by poor ones. ME and SVDDPs can explore and find better solutions with the exploration mechanisms. To examine the effect of temperature $\alpha$ and the number of modes (trajectories), we performed the same experiments, with different values of them. The results of the experiments are provided in Table \ref{tab:TO_alpha_mode}. For SVDDP, we also changed step size $\epsilon$s. We plot the evolution of the mean cost over iteration in Fig. \ref{fig:TO_convergence} to compare the convergence speed of DDPs.\\
\underline{2D Car:}
The state $x\in \Rb^{3}$ is composed of positions and orientation, and the control $u\in \Rb^{2}$ is composed of transitional and angular velocity. All algorithms can find the global solution when the number of trajectories increases regardless of $\alpha$. This is because the state and control space is explored more as the number of modes increases. The performance of SVDDP is not as good as MEDDPs with low $\alpha$ ($\leq 20$) and small $\epsilon$ ($=1$). However, with sufficiently large values, SVDDP improves its performance and surpasses MEDDPs. The reason for the result is that these values amplify the repulsive force for exploration as in \eqref{eq:Hss_SVDDP} and \eqref{eq:SV_ddp_update}. For convergence speed, the SVDDP outperforms the MEDDPs.\\
\underline{Quadrotor:}
This system has positions, angles, and their derivatives as a state, and the control is force generated by four rotors. Hence, it has $x\in\Rb^{12}, u\in\Rb^{4}$. We can observe two similar tendencies to the previous example. All algorithms improve performance as the number of modes increases and SVDDP starts outperforming MEDDPs in high temperatures and enough step size. The convergence speed of SVDDP is slightly better, but not as significant as in the case of the 2D car.\\
\underline{Arm:}
An arm has 7 joints. Adding the position of the end effector, it has $x\in \Rb^{17}, u\in\Rb^{7}$.
In addition to the previous observations, we see that a pair of too high $\alpha$ and $\epsilon$ does not help to improve the performance of SVDDP. This indicates that a strong repulsive force can significantly disturb all the optimization efforts.

\textbf{Model Predictive Control:} \label{MPC_experiments}
We perform experiments in MPC with the 2D car and quadrotor system used in TO for similar reaching tasks. With 2D car and quadrotor, we compare seven algorithms of MPC DDPs and one of the most successful sampling-based algorithms MPPI to examine the performance of our algorithms. We have DDP, UG-MEDDP, MG-MEDDP, and SVDDP in MPC mode for DDP-based algorithms. For MPPI, we have UG-, MG-, and SV-MPPI. Details of the experiments are provided in Appendix \ref{sec:app_MPC_results}. We use \cite{nni2021} to tune MPPIs parameters such as covariance and temperature, etc. 

In the experiments, a trajectory is considered successful when it reaches a ball centered at the target with radius 0.5, and does not hit the obstacles. The experiments are performed ten times per algorithm and environment, except for the deterministic normal DDP. This means that for each algorithm we will have 100 runs from which we collect statistics on the performance.  To show that DDP is the effective method for ME second-order optimization, we also perform the same experiment by discarding the feedback gain $K$ in \eqref{eq:delta-u-star-main} and \eqref{eq:SV_ddp_update} when sampling trajectories in ME and SVDDPs. This is because providing feedback gain is a unique property of DDP in contrast to other second-order methods. 

\begin{figure}[!h]
\centering
    \begin{subfigure}[h!]{0.48\linewidth} 
    \includegraphics[trim={0.1cm 0.0cm 6.5cm 0.0cm},clip,width=\textwidth]{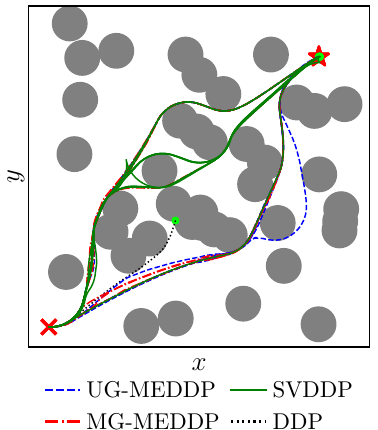}
    \caption{2D Car MPC DDPs}
    \label{fig:car_mpc_ddp}
    \end{subfigure}
    \begin{subfigure}[h!]{0.48\linewidth} 
    \includegraphics[trim={6.5cm 0.0cm 0cm 0.0cm},clip,width=\textwidth]{fig/MPC/car/Paper/Twod_car_DDP_MPPI_Grid10.pdf}
    \caption{2D Car MPPIs}
    \label{fig:car_mppi}
    \end{subfigure}
    \begin{subfigure}[h!]{0.495\linewidth} 
    \includegraphics[trim={1.3cm 0.1cm 1.3cm 1.3cm},clip,width=\textwidth]{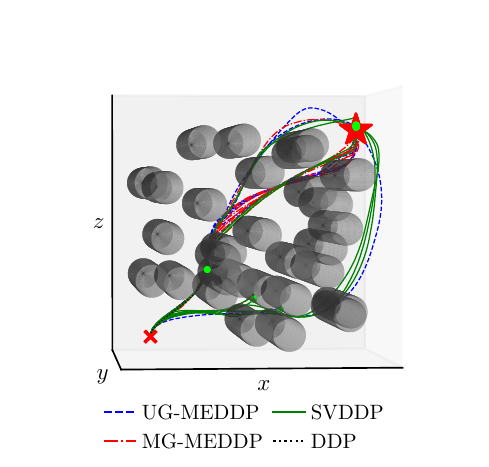}
    \caption{Quadrotor MPC DDPs}
    \label{fig:quad_mpc_ddp}
    \end{subfigure}
    \begin{subfigure}[h!]{0.495\linewidth} 
    \includegraphics[trim={1.3cm 0.1cm 1.3cm 1.3cm},clip,width=\textwidth]{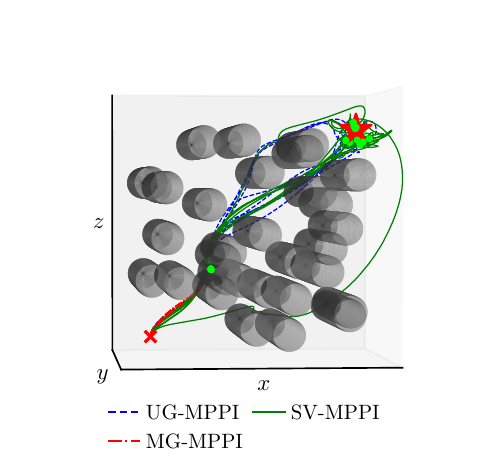}
    \caption{Quadrotor MPPIs}
    \label{fig:quad_mppi}
    \end{subfigure}
\caption{Comparison of MPPIs and MPC DDPs for 2D Car and Quadrotor in one of the ten obstacle fields. {\color{red}$\bm{\times}$}:Initial positions. {\color{red}$\bigstar$}:Target positions. {{\color{lime}{\CircleSolid}}}: position at terminal time step.}
\label{fig:cluttered_env_mpc_ddps}
\end{figure}

\underline{2D Car:}
One of the obstacle fields and results are shown in Fig. \ref{fig:car_mpc_ddp} and Fig. \ref{fig:car_mppi}. The initial point and the target are [0, 0] and [5, 5]. The normal DDP get stack while MEDDPs and SVDDPs are able to hit the target. MPPIs can explore the fields, but not as good as ME and SV DDPs and they tend to require longer time even when they escape local minima.
Statistics on the performance of the algorithms with respect to the success rate are provided in Table \ref{tab:forest_twod_car}. Here, the average constraint violation is computed over trajectories that hit the target with constraint violation. The constraint is computed by $({\text{center of obstacle}-\text{position of vehicle}})^{2} - ({\text{radius of obstacle}})^{2}.$ The results illustrate that SVDDP outperforms the other methods in all metrics. In addition, the results indicate that a proper choice of a margin of constraints increases the success rate of DDPs even further. Without feedback gains in sampling, the algorithms reduce the success rate as SVDDP: 74 \%, MG-MEDDP: 71 \% and  UG-MEDDP 68 \%, which highlights the importance of feedback gains when sampling trajectories. The statistics of the computational time of one MPC call over ten different runs are presented in Table \ref{tab:cmp_time_car}. The algorithms are implemented with JAX \cite{jax2018github} and executed with Nvidia GeForce RTX 3080 GPU. Here, the prediction horizon amounts to 1.2 sec as in section \ref{sec:app_MPC_2D_car}. Although SVDDP takes longer time than other methods, it is still short enough for real-time implementation. Finally, we emphasize that DDPs can generate smoother control trajectories as presented in Fig. \ref{fig:app_cntl_smoothness}, which shows one of the two control commands with control bounds. \\
\underline{Quadrotor:}
One of the obstacle fields and results are shown in Fig. \ref{fig:quad_mpc_ddp} and Fig. \ref{fig:quad_mppi}. The obstacles are cylinders whose axes aligns with the $y$ axis. The initial and target points are [0, 0, 0] and [5, 0, 5]. The statistics of the results and computational time are provided in Table \ref{tab:forest_quad}, and \ref{tab:cmp_time_quad}. Again, SVDDP outperforms other methods in success rate. It takes longer for computation, but still fast enough for real-time application. Without feedback gains in sampling,  we obtained the following success rate. SVDDP: 45 \%, MG-MEDDP: 56 \% and  UG-MEDDP 52 \%, which indicates the importance of the feedback gain even more in an unstable system.


\begin{table}[h]
\caption{Statistics for 2D car and quadrotor navigation task in randomized obstacle environments.}
\label{tab:temps}
\footnotesize
\begin{subtable}[h]{0.48\textwidth}
\centering
\caption{2D Car}
\label{tab:forest_twod_car}
\resizebox{\columnwidth}{!}{
\begin{tabular}{c|c|c|c}
Method   & \makecell{Success\\ rate [\%]} & 
\makecell{Success\\ with const.\\violation[\%]} & 
\makecell{avg. const. \\violation}
\\ \hhline{=|=|=|=}
SVDDP    & 84  & 100 &    \num{9.49e-5}      \\
MG-MEDDP & 72  & 93  &    \num{1.91e-4}      \\
UG-MEDDP & 79  & 98  &    \num{3.36e-4}      \\
DDP      & 10  & -   &  -       \\ \hline
SVMPPI   & 39  & -   &  -        \\
MG-MPPI  & 8   & -   &  -        \\
UG-MPPI  & 6   & -   &  -       
\end{tabular}}
\end{subtable}
\begin{subtable}[h]{0.48\textwidth}
\caption{Quadrotor}
\label{tab:forest_quad}
\resizebox{\columnwidth}{!}{
\begin{tabular}{c|c|c|c}
Method   & \makecell{Success\\ rate [\%]} & 
\makecell{Success \\ with const.\\ violation[\%]} & 
\makecell{avg. const. \\violation}
\\ \hhline{=|=|=|=}
SVDDP    & 92  & 94 &    \num{6.94e-6}      \\
MG-MEDDP & 71  & 90  &   \num{1.02e-5}       \\
UG-MEDDP & 70  & 86  &   \num{2.45e-5}       \\
DDP      & 30  & -   &  -       \\ \hline
SVMPPI   & 69  & -   &  -        \\
MG-MPPI  &  9  & -   &  -        \\
UG-MPPI  & 72  & -   &  -       
\end{tabular}}
\end{subtable}
\end{table}

\begin{table}[h]
\caption{Computation times for one MPC call averaged over 10 different runs.}
\label{tab:cmp_time}
\begin{subtable}[h]{0.45\textwidth}
\footnotesize
\caption{2D car}
\label{tab:cmp_time_car}
\centering
\begin{tabular}{c|cc}
Method   & Mean[s] & Std.[s] \\ \hline
SVDDP    & \num{4.02e-2}  & \num{9.37e-3} \\
MG-MEDDP & \num{3.43e-2} & \num{1.19e-3} \\
UG-MEDDP & \num{3.45e-2} & \num{8.57e-4}
\end{tabular}
\end{subtable}
\hfill
\begin{subtable}[h]{0.45\textwidth}
\footnotesize
\caption{Quadrotor}
\label{tab:cmp_time_quad}
\centering
\begin{tabular}{c|cc}
Method   & Mean[s] & Std.[s] \\ \hline
SVDDP    & \num{7.36e-2}  & \num{4.15e-2} \\
MG-MEDDP & \num{4.11e-2} & \num{1.66e-2} \\
UG-MEDDP & \num{4.76e-2} & \num{2.54e-2}
\end{tabular}
\end{subtable}
\end{table}

\textbf{Limitations:} To find a better local solution by exploration, especially in MPC, a sufficiently long prediction time horizon is required. The choice of the time horizon can depend on the environment, such as the density of obstacles, and dynamics of the system, including control limits. Therefore, one may need to adjust the time horizon to enjoy the exploration capability of SV and MEDDPs. The temperature $\alpha$ may also need to be tuned. We experimentally verified that SVDDP works stably in a wide range of $\alpha$, e.g., $\alpha \in [10,40]$ for 2D car.  However, the range depends on the environment and the system. Therefore, one may need to tune $\alpha$ to enjoy the performance of SVDDP.

\section{Conclusion}\label{sec:conclusion}
In this paper, we have derived SVDDP based on the MEDDP formulation and SVNM, which has superior exploration capability to the original MEDDPs. Moreover, we incorporate state and control constraints with the barrier method to achieve constrained optimization with SVDDP. We compare SVDDP with a constrained version of MEDDP in TO mode in three different systems, and also with MPPIs in two different systems in MPC mode. The results suggest that SVDDP can outperform MPPI in finding optimal solutions and in generating smooth state and control trajectories. Furthermore, SVDDP shows better exploration capability compared to MEDDPs. Future work includes hardware implementation with multimodal cost to validate the algorithm. Although we used the barrier function to represent constraints in this work, alternative optimization techniques can also be incorporated into SVDDP. In addition, while SVDDP in its current form requires differentiability, extensions that handle non-smooth dynamics \cite{LELIDEC2024101468}  and non-smooth cost functions \cite{DDP_NonSmoothObj} are straightforward. Finally, theoretical analysis of the disruptive mechanism in \textit{Messy} Dynamic Optimization and SVDDP is of great importance.  

\clearpage


	



\clearpage
\acknowledgments{If a paper is accepted, the final camera-ready version will (and probably should) include acknowledgments. All acknowledgments go at the end of the paper, including thanks to reviewers who gave useful comments, to colleagues who contributed to the ideas, and to funding agencies and corporate sponsors that provided financial support.}

\appendix
\renewcommand{\theequation}{\thesection.\arabic{equation}}
\renewcommand{\thefigure}{\thesection.\arabic{figure}} 
\renewcommand{\thetable}{\thesection.\arabic{table}}
\setcounter{figure}{0} 
\setcounter{table}{0} 
\setcounter{equation}{0}
\section{Maximum Entropy Differential Dynamic Programming}\label{sec:app_DDP}
This section provides derivation of Differential Dynamic Programming (DDP) \citep{Jacobson1970ddp} and Maximum Entropy DDP (MEDDP) \citep{so2022maximum} for completeness.
\subsection{Optimal Stochastic Policy of Entropic Regularized Optimal Control Problem}

To solve the problem above, we construct a Lagrangean
\begin{align}
\mathcal{L} = \int \Pi(U) [J + \alpha \ln(\Pi(U)]{\mathrm{d}}U + \lambda (1- \int \Pi(U){\mathrm{d}}U),
\end{align}
where $\lambda$ is a Lagrangian multiplier for the constraint for $\Pi$ being a valid probability distribution. 
We take a functional derivative with respect to $\Pi$ and set it zero, obtaining
\begin{align}\label{eq:pi_general_ent_opt}
\notag
& \  J + \alpha \ln (\Pi(U)) + \alpha -\lambda = 0  \\
\Leftrightarrow &  \ \Pi(U) = \exp[(-J + \lambda-\alpha)/{\alpha}]
= \exp[- J/\alpha]\exp[{\lambda}/{\alpha}-1].
\end{align}
We eliminate $\lambda$ from the equation above using the constraint for $\Pi$ as follows:
\begin{align*}
\int\Pi(U) \mathrm{d}U = \int\exp[-{J}/{\alpha}] \mathrm{d}U \exp[{\lambda}/{\alpha}-1] = 1  \\ 
\Leftrightarrow  \exp[{\lambda}/{\alpha}-1] = \Big[\int\exp[-{J}/{\alpha}] \mathrm{d}U\Bigg]^{-1}.
\end{align*}
By plugging this into \eqref{eq:pi_general_ent_opt}, we obtaing
\begin{align}
\Pi(U) = \frac{1}{Z}\exp\Big[-\frac{J}{\alpha}\Big], \quad \text{with} \quad Z = \int\exp\Big[-\frac{J}{\alpha}\Big] \mathrm{d}U
\end{align}

\subsection{Differential Dynamic Programming}
Consider a trajectory optimization problem for a dynamical system with state $x\in \mathbb{R}^{n_{x}}$ and control $u\in \mathbb{R}^{n_{u}}$. 
Given the deterministic dynamics
$x_{t+1} = f(x_{t}, u_{t})$,
The trajectory optimization problem is formulated as follows:   
\begin{equation}\label{eq:ocp_cost_app}
J(X,U) = \Phi(x_{T}) +  \sum_{t=0}^{T-1}l_{t}(x_{t},u_{t}) , \quad \text{s.t.}  \quad x_{t+1} = f(x_{t}, u_{t})
\end{equation}
where $X = [x_{0}, \cdots, x_{T}], U = [u_{0}, \cdots, u_{T-1}]$. The scalar-valued functions $l_{t}(\cdot,\cdot)$, $\Phi(\cdot)$, and $J(\cdot, \cdot)$ denote the running, terminal, and total cost of the problem, respectively. The cost-to-go at time step $t=i$, i.e., cost starting from $t=i$ to $N$ denoted by ${J}_{i}(X_{i},U_{i})$ is given by
\begin{equation*}
     {J}_{i}(X_{i},U_{i}) = \big[\sum_{t=i}^{N-1} l_{t}(x_t,u_t)\big] + \Phi(x_{N})
\end{equation*}
with $X_{i} =[{x}_{i},\dots,{x}_{N}]$, $U_{i} =[{u}_{i},\dots,{u}_{N-1}]$. The value function is defined as the minimum cost-to-go in each state and time step $t$ via $V_t({x}_t):=\min_{{u}_t}J_{t}({X_{t}},{U_{t}})$.
Given Bellman's principle of optimality that provides the following rule
\begin{equation}\label{eq:bellman}
    V_t(x_t) =  \min_{{u}_t} [\underbrace{l_{t}({x}_t,{u}_t) + V_{t+1}({x}_{t+1})}_{Q_{t}(x_{t},u_{t})}], 
\end{equation}
DDP finds local solutions to \eqref{eq:ocp_cost_app} by expanding both sides of   \eqref{eq:bellman} about nominal trajectories $\bar{X}$ and $\bar{U}$. 
To proceed, we first take quadratic expansions of $Q_t$ about $\bar{x}_{t}$ and $\bar{u}_{t}$ with the deviation $\delta x_t=x_t -\bar{x}_{t}$, $\delta u_t=u_{t} - \bar{u}_{t}$, obtaining
\begin{equation}\label{eq:Q_quad_expand}
 Q(x_{t},u_{t}) \approx  Q(\bar{x}_{t}, \bar{u}_{t}) + Q_{u}^{\tr}\delta u_{t} + Q_{x}^{\tr}\delta x_{t} +\frac{1}{2}\delta u^{\tr}Q_{uu}\delta u_{t} + 
\delta u^{\tr}Q_{ux}\delta x_{t} 
 +\frac{1}{2}\delta x_{t}^{\tr}Q_{xx}\delta x_{t},
\end{equation}
where we drop $t$ for $Q$.
Then, in the definition of the $Q$ given in \eqref{eq:bellman}, we perform a similar expansion and get the following expression:
\begin{align}\label{eq:Q_quad_expand_lV}
l(x_{t},u_{t}) + V(x_{t+1}) &\approx l(\bar{x}_{t}, \bar{u}_{t}) +  V(\bar{x}_{t+1}) + l_{u}^{\tr}\delta u_{t} + l_{x}^{\tr}\delta x_{t} +\frac{1}{2}\delta u_{t}^{\tr}l_{uu}\delta u_{t} + 
\delta u_{t}^{\tr}l_{ux}\delta x_{t} \\\notag
&\quad +\frac{1}{2}\delta x_{t}^{\tr}l_{xx}\delta x_{t} 
+ V_{x}(x_{t+1})^{\tr} \delta x_{t+1} + \frac{1}{2} \delta x_{t+1} ^{\tr}V_{x}(x_{t+1}) \delta x_{t+1}.
\end{align}
Here, we drop $t$s from $l(\cdot, \cdot)$ and $V(\cdot)$, but they can be recovered from those of the arguments.
From the above equation, we eliminate $\delta x_{t+1}$ by plugging in the linearized dynamics given by
\begin{align}\label{eq:linearized_dyn}
\delta x_{t+1} =  f_{x} \delta x_{t} + f_{u} \delta u_{t},
\end{align}
where $f_{x}$ and $f_{u}$ are Jacobian matrices for state and control, respectively.
Mapping the terms on both sides of $Q$ (\eqref{eq:Q_quad_expand} and \eqref{eq:Q_quad_expand_lV}), we obtain the derivatives of $Q$ (evaluated on $\bar{X}$ and $\bar{U}$) as follows. 
\begin{align}\label{eq:Qexpanded_derivative}
    {Q}_{{u},t}&={l}_{{u}}+{f_u}^{\tr}{V}_{{x},t+1}, \quad {Q}_{{x},t}={l}_{{x}}+{f_x}^{\tr}{V}_{{x},k+1} \\\notag
    {Q}_{{uu},t} &={l}_{{uu}}+{f_u}^{\tr}{V}_{{xx},t+1}{f_u}, \quad
    {Q}_{{xu},t}  ={l_{{xu}}}+{f_x}^{\tr}{{V}_{{xx},t+1}}{f_u} \\\notag 
    {Q}_{{xx},t}  &={l}_{{xx}}+{f_x}^{\tr}{V}_{{xx},t+1}{f_x}.
\end{align}
Assuming that $Q_{uu}$ is Positive Definite (PD), we can explicitly optimize $Q$ approximated in \eqref{eq:Q_quad_expand} with respect to $\delta{u}_{t}$ by computing a partial derivative of \eqref{eq:Q_quad_expand} with respect to $\delta u_{t}$ and setting it zero. This minimization yields the following local optimal control law
\begin{align}\label{eq:delta-u-star}
    & \delta {u}^{\ast}_t={\kappa}_{t}+{K}_{t}\delta {x}_{t}, \quad 
    \text{with}\quad \kappa_{t}= -{Q}^{-1}_{{uu,t}}{Q_{{u,t}}},\ {K}_{t} = -{Q}^{-1}_{{uu,t}}{Q_{{ux,t}}},
\end{align}
where $\kappa$ and $K$ are known as feedforward and feedback gains, respectively. Now, $\delta {u}_{t}^{\ast}$ is computed using $V_{x,t+1}$ and $V_{xx,t+1}$, which are in one time step ahead. To propagate derivatives of $V_{t}(x)$ back in time, we consider the quadratic expansion of $V_{t}(x)$, that is, 
\begin{equation}\label{eq:V_quad_expand}
V(x_{t}) = V(\bar{x}_{t}) + V_{x}(\bar{x}_{t})^{\tr}\delta x_{t} + \frac{1}{2}\delta x_{t}^{\tr}V_{xx}(\bar{x}_{t})\delta x_{t}, 
\end{equation}
and equate the equation with the quadratic expansion of $Q$ in \eqref{eq:Q_quad_expand} through \eqref{eq:bellman}. Since we now have a solution of the $\min_{u_{t}}$ in the right-hand side of \eqref{eq:bellman}, by substituting $\delta u^{\ast}$ for $\delta u$, the $\min$ operator vanishes. This allows us to compare the coefficients of $\delta x_{t}$, giving the backward recursions:
\begin{align} \label{eq:value_riccati}
    V_{{x},t}={Q_{{x},t}}-{Q_{{xu},t}}{Q_{{uu},t}^{-1}}{Q_{{u},t}},\quad
    V_{{xx},t}={Q_{{xx},t}}-{Q_{{xu},t}}Q_{{uu},t}^{-1}{Q_{{ux},t}},
\end{align}
with the terminal condition
\begin{equation}\label{eq:value_terminal_cnd}
V_{N}(x_{T}) = \Phi(x_{T}), \quad V_{x,T} = \Phi_{x}(x_{T}), \quad V_{xx,T} = \Phi_{xx}(x_{T}). 
\end{equation}

DDP has backward and forward passes. In the backward pass, the derivatives of $Q$ functions in \eqref{eq:Qexpanded_derivative}, gains in \eqref{eq:delta-u-star}, and derivatives of the value function in \eqref{eq:value_riccati} are computed backward in time. In the forward pass, the new control $\bar{u}_{t} + \delta u^{\ast}_{t}$ is propagated forward in time to give a new pair of nominal trajectories, typically with a backtracking line search.

\subsection{Maximum Entropy DDP}
This section provides a derivation of Maximum Entropy DDP (MEDDP). Starting with the entropic regularized cost and corresponding optimal policy, we derive unimodal and multimodal policies of MEDDPs. Finally, we provide a schematic of the algorithms.

\subsubsection{Optimal Policy}
From the entropic regularized cost \eqref{eq:ME_DDP_objective} in the main paper, we obtain the following optimization problem to compute the optimal polity $\pi^{\ast}$:
\begin{align}\label{eq:sup_max_entropy_with_cnst}
\min_{\pi} \mathbb{E}_{u \sim \pi(\cdot|x)}[l(u,x) + V(x')] - \alpha H[\pi(\cdot | x)],  \quad \text{subject to} 
\int \pi(u|x)\mathrm{d}u = 1,
\end{align}
where we dropped time instance $t$ and denote $x_{t+1}$ as $x'$ for readability. $\pi(\cdot | x)$ means that it originally was a function of $u$, but $u$ vanishes after expectation was taken over $u$. 
To solve this equation, we first form a Lagrangian $\mathcal{L}$ with Lagrangian multiplier $\lambda \in \mathbb{R}$ as 
\begin{align*}
\mathcal{L} &= \mathbb{E}_{u \sim \pi}[l(u,x) + V(x')] - \alpha H[\pi(\cdot | x)] + \lambda 
 \Big(1-\int \pi(u|x)\mathrm{d}u \Big),\\
&= \mathbb{E}_{u \sim \pi}[l(u,x) + V(x') + \alpha \ln \pi(u | x) - \lambda] + \lambda \\
&= \int_{\mathbb{R}^{n_{u}}}\pi(u|x) [Q(x,u) + \alpha \ln \pi(u|x) - \lambda]\mathrm{d}u + \lambda,
\end{align*}
where we used the definition of $Q$ in \eqref{eq:bellman} in the last equality. 
From the first-order optimality condition, taking functional derivative of $\mathcal{L}$ with respect to $\pi$ and setting it zero, we have 
\begin{align*}
Q(x,u) + \alpha \ln \pi(u|x) - \lambda + \alpha = 0 \ 
\Leftrightarrow \ \ln \pi(u|x) = -\frac{1}{\alpha}[Q(x,u)-\lambda] -1.
\end{align*}
Thus, the optimal policy $\pi^{\ast}$ is obtained as
\begin{align}\label{eq:MEDDP_pi_star_with_lam}
\pi^{\ast}(u|x) = \exp \Big[ -\frac{1}{\alpha}(Q(x,u) - 
\lambda) - 1\Big].
\end{align}
From the constraint in \eqref{eq:sup_max_entropy_with_cnst}, we eliminate $\lambda$ from the above equation. First, we plug \eqref{eq:MEDDP_pi_star_with_lam} into the constraint, having
\begin{align*}
1 = \int \pi^{\ast}(u|x) \mathrm{d}u =
\int \exp \Big[ -\frac{1}{\alpha} Q(x,u) \Big]
\exp\Big[ 
\frac{\lambda}{\alpha}-1
\Big] \mathrm{d}u \\
\Leftrightarrow
\exp \Big[ \frac{\lambda}{\alpha} -1 \Big] = \Big( \int 
\exp\Big[ \frac{1}{\alpha}Q(x,u) \mathrm{d}u\Big]\Big)^{-1} = Z^{-1},
\end{align*}
where the partition function is given by
\begin{align}\label{eq:MEDDP_Z}
Z = \int \exp\Big[ 
-\frac{1}{\alpha}Q(x,u)
\Big] \mathrm{d}u = \int \exp\Big[ 
-\frac{1}{\alpha}(l(x,u) + V(x'))
\Big]\mathrm{d}u.
\end{align}
By plugging this into  \eqref{eq:MEDDP_pi_star_with_lam}, the optimal policy is obtained as 
\begin{align}\label{eq:sup_pi_star}
\pi^{\ast}(u|x) = \exp \Big[ -\frac{1}{\alpha}Q(x,u)\Big]
\exp \Big[ \frac{\lambda}{\alpha}-1 \Big] = \frac{1}{Z}\exp\Big[ -\frac{1}{\alpha}Q(x,u)\Big]. 
\end{align}

To derive the relationship between the value and partition functions in 
\eqref{eq:me_ddp_pi_star_uni}, we transform the equation above to get
\begin{align*}
Z\pi^{\ast}(u|x) = \exp \Big[-\frac{1}{\alpha} Q(x,u)\Big]
\Leftrightarrow \ 
\ln[Z \pi^{\star}(u|x)] = - \frac{1}{\alpha}Q(x,u) \\
\Leftrightarrow \ Q(x,u) = -\alpha \ln(Z \pi^{\ast}(u|x)) = -\alpha (\ln Z + \ln \pi^{\ast}(u|x)).
\end{align*}

In \eqref{eq:sup_max_entropy_with_cnst}, since $H[\pi(\cdot|x)]$ is a constant and not a function of $u$, the term can be inside of the expectation, which gives the entropic regularized version of value function:
\begin{align*}
V(x)
&= \min_{\pi}\big[ \mathbb{E}_{u\sim\pi}[l(x, u) + V(x')] - \alpha \ln H[\pi(\cdot|x)] \big].
\end{align*}
From the above two equations, we get the following relationship.
\begin{align}\label{eq:sup_V_islog_Z}
V(x) &= \mathbb{E}_{\pi^{\ast}}[Q-\alpha H[\pi^{\ast}]] \\\notag
& = \int {\pi^{\ast}(u|x)} [-\alpha \ln{Z} - \alpha \ln\pi^{\ast}(u|x)+\alpha \ln \pi^{\ast}(u|x)] \mathrm{d}u \\\notag
&= -\alpha \ln Z \int \pi^{\ast}(u|x)\mathrm{d}u \\\notag
&= -\alpha \ln Z(x).
\end{align}
 
\subsubsection{Unimodal policy}
This section provides the derivation of the unimodal Gaussian policy of MEDDP.
Consider a quadratic approximation of $Q$ around a pair of nominal trajectories as in \eqref{eq:Q_quad_expand}. We write its nominal terms with the running cost and value function, and its deviation with $Q$. In addition, we complete the square of the deviation term with respect to $\delta u$:
\begin{align}\label{eq:sup_MEDDP_expanded_Q}
Q(x,u) &= l(\bar{x}, \bar{u}) + V(\bar{x}') + \delta Q(\delta x, \delta u),\\\notag
\text{with} \ \
\delta Q &= Q_{x}^{\tr}\delta x + Q_{u}^{\tr}\delta u + \frac{1}{2}\delta x ^{\tr}Q_{xx}\delta x+ \delta x ^{\tr}Q_{xu}\delta u + \frac{1}{2}\delta u ^{\tr}Q_{uu}\delta u\\\notag
&= \frac{1}{2} \delta u^{\tr}Q_{uu}\delta u + \big[Q_{u} + Q_{ux}\delta x  \big]^{\tr}\delta u + Q_{x}^{\tr}\delta x
+ \frac{1}{2}\delta x ^{\tr}Q_{xx}\delta x\\\notag
&= \frac{1}{2}(\delta u - \delta u^{\ast})^{\tr}Q_{uu}(\delta u - \delta u^{\ast}) \underbrace{-\frac{1}{2} \delta u^{\ast \tr}Q_{uu}\delta u^{\ast} + Q_{x}^{\tr}\delta x
+ \frac{1}{2}\delta x ^{\tr}Q_{xx}\delta x}_{\text{not a function of $u$,}}
\end{align} 
as indicated in the equation above, in $\delta Q$, only the first term is a function of $u$.
$\delta u^{\ast}$ is given in \eqref{eq:delta-u-star}, but it is just a constant to complete the square, but not a minimizer. Plugging \eqref{eq:sup_MEDDP_expanded_Q} into the optimal policy in  \eqref{eq:sup_pi_star}, we obtain
\begin{align*}
\pi^{\ast}(u|x) \propto \exp \Big[ - \frac{1}{\alpha} Q(x,u)\Big] \propto \exp \Big[- {\frac{1}{2\alpha}}
(\delta u - \delta u^{\star})^{\tr}Q_{uu}(\delta u - \delta u^{\star})\big] \\
\propto
\exp \big[-{\frac{1}{2}}
(\delta u - \delta u^{\ast})^{\tr}\Sigma^{-1}(\delta u - \delta u^{\ast})\big], \quad \text{with} \quad \Sigma = \alpha Q_{uu}^{-1},
\end{align*}
which indicates that the optimal policy is normally distributed with mean $\delta u^{\ast}$ and covariance $\alpha Q_{uu}^{-1}$. We see that when temperature $\alpha$ is high, the optimal policy has a higher variance, and thus becomes noisy.

\subsubsection{Multimodal policy}
This section provides the derivation of MEDDP's multimodal Gaussian policy, which is given as a weighted sum of unimodal Gaussian policies. This distribution is also known as a mixture of Gaussians, and the set of weights is a categorical probability mass function that determines which mode (unimodal Gaussian distribution) is used when sampling data from the distribution. In multimodal policy, the weights are proportional to the value function of the trajectories associated with unimodal policies. In this section, we first compute the value function and partition function of a unimodal policy associated with a single trajectory. Then, we show that the transformation of these functions and the policies have additive structures. Finally, we show the form of the weight in \eqref{eq:sup_weight} and explain how to compute it. 

\textbf{Value function and partition function:}
To compute the value function,
we use the same quadratic approximation of $Q$ in \eqref{eq:sup_MEDDP_expanded_Q} and group the terms that are not a function of $\delta u$ and name it $\tilde{V}(x)$:
\begin{align}\label{eq:sup_Q_expand}
\frac{1}{\alpha}Q(x,u) &= \frac{1}{\alpha}
[\underbrace{V(\bar{x}') + l(\bar{x}, \bar{u})-\frac{1}{2} \delta u^{\ast \tr}Q_{uu}\delta u^{\ast} +  Q_{x}^{\tr}\delta x
+ \frac{1}{2}\delta x ^{\tr}Q_{xx}\delta x}_{\tilde{V}(x)} \\\notag
& \quad \quad +\frac{1}{2}(\delta u - \delta u^{\ast})^{\tr}Q_{uu}(\delta u - \delta u^{\ast})].
\end{align}
The name of the term $\tilde{V}(x)$ comes from the fact that it is obtained by substituting $\delta u = \delta u ^{\ast}$ into \eqref{eq:sup_Q_expand}, which is a value function of normal DDP for the state $x$ (at time $t$).
This is verified by the following equation. 
\begin{align*}
\text{Value function of DDP for }x_{t} = \min_{u_{t}}Q(x_{t},u_{t}) = Q(x_{t},u_{t}^{\ast}) = Q(x_{t},\bar{u}_{t} + \delta u_{t}^{\ast}) = \tilde{V}(x_{t}) .
\end{align*}
The quadratic approximation also allows us to compute $Z$ in \eqref{eq:MEDDP_Z} using Gaussian integral, that is, an integral has a form of 
\begin{align*}
\int \exp \Big[ \frac{1}{2}y^{\tr}Ay\Big] \mathrm{d}y = [(2\pi)^{n_{y}}|A^{-1}|]^{1/2}, 
\end{align*}
with $y\in \Rb^{n_{y}}$, and a symmetric PD matrix $A\in\Rb^{n_{y} \times n_{y}}$. Applying this integral to \eqref{eq:MEDDP_Z} with \eqref{eq:sup_Q_expand}, we have
\begin{align*}
Z &= \int \exp\Big[ -\frac{\tilde{V}(x)}{\alpha}\Big]  \exp\Big[
-\frac{1}{2}(\delta u - \delta u^{\ast})^{\tr}\frac{Q_{uu}}{\alpha}(\delta u - \delta u^{\ast})
\Big] \mathrm{d}u \\
&= \exp\Big[ -\frac{\tilde{V}(x)}{\alpha}\Big][(2\pi)^{n_{u}}|\alpha Q_{uu}^{-1}|]^{1/2}.
\end{align*}
Substituting this $Z$ back into \eqref{eq:sup_V_islog_Z}, the value function is obtained as the sum of value function of normal DDP and a term from the Gaussian policy which we call $V_{H}$:
\begin{align}\label{eq:me_ddp_value_function_lnz}
V(x) &= -\alpha \ln Z \\\notag
&=  -\alpha \big[ 
-\frac{\tilde{V}(x)}{\alpha} + \frac{1}{2}\ln((2\pi)^{n_{u}}|\alpha Q_{uu}^{-1})) \big]\\\notag
&= \tilde{V}(x) + \underbrace{\frac{\alpha}{2}\big[
\ln |Q_{uu}| -n_{u} \ln(2\pi \alpha)
\big]}_{V_{H}}.
\end{align}

\textbf{Additive structure:} Instead of performing a single quadratic approximation of the cost around a pair of nominal trajectories, we consider the combined approximation. We start from the terminal cost $\Phi(x_{T})$ with $N$ sets of trajectories:
$
X^{(n)} = [x_
{0}^{(n)}, \cdots, x_{T}^{(n)}], \quad 
U^{(n)} = [u_
{0}^{(n)}, \cdots, u_{T-1}^{(n)}],
$
with $n=1 ,\cdots, N$. Using the LogSumExp operation and the terminal condition \eqref{eq:value_terminal_cnd}, we first approximate the value function at terminal timestep $T$ using $N$ trajectories as 
\begin{align}\label{eq:MEDDP_combined_terminal}
V(x_{T}) = \tilde{\Phi}(x_{T})
= -\alpha \ln \sum_{n=1}^{N} \exp\Big[
-\frac{1}{\alpha}\Phi^{(n)}(x_{T})
\Big],
\end{align}
where $\Phi^{(n)}(x)$ denotes $\Phi(\cdot)$ for $x^{(n)}$. To proceed, we introduce the exponential transformation
$z_{t} = \mathcal{E}_{\alpha }[V(x_{t})]$, with $\mathcal{E}_{\alpha}(y) = \exp[
-y/{\alpha}]$, which transforms the cost and value functions into reward and desirability functions, respectively. 
\begin{align*}
r_{t} = \mathcal{E}_{\alpha }[l_{t}], \quad  
r_{T} = \mathcal{E}_{\alpha }[\Phi(x_{T})], \quad  
z_{t} = \mathcal{E}_{\alpha }[V(x_{t})].
\end{align*}
From \eqref{eq:MEDDP_combined_terminal} and \eqref{eq:value_terminal_cnd},
at the terminal time step, $z$ (the approximation by $z^{(n)}$s) is written as a sum of $z^{(n)}$s due to the property of exponential function as below.
\begin{align}
z_{T}(x_{T}) &= 
\exp \Big[
-\frac{1}{\alpha}\Big(-\alpha \ln \sum_{n=1}^{N}\exp \Big[ -\frac{1}{\alpha}\Phi^{(n)}(x)\Big]\Big)
\Big]\\\notag 
&=\sum_{n=1}^{N} \exp\Big[
-\frac{1}{\alpha}\Phi^{(n)}(x) 
\Big]
= \sum_{n=1}^{N}z_{T}^{(n)}(x_{T}).
\end{align}
This additive structure holds not only at the terminal time step, which can be proved by the following induction. We use the following two relationships for the induction.
\begin{align}\label{eq:sup_z_equal_Z}
z(x) = \mathcal{E}_{\alpha}[V(x)] = \exp \Big[
-\frac{1}{\alpha}(-\alpha \ln Z(x))
\Big] = Z(x),
\end{align}
and 
\begin{align}\label{eq:sup_Z_rz}
Z(x) 
&=\int \exp \Big[
-\frac{1}{\alpha}(l(x,u)+ V(x'))\Big] \mathrm{d}u\\\notag
&=\int \exp \Big[
-\frac{1}{\alpha}l(x,u)\Big]
\exp \Big[ 
-\frac{1}{\alpha}V(x')
\Big] \mathrm{d}u\\\notag
&= \int r(x,u) z(x')\mathrm{d}u.
\end{align}
Assume the additive structure hods at time instance $t+1$, that is
\begin{align*}
z(x') = \sum_{n=1}^{N}z^{(n)}(x'). 
\end{align*}
Substitution of the equation above into \eqref{eq:sup_Z_rz}, at time $t$, we get 
\begin{align}\label{eq:sup_additive_z}
z(x) = \int \sum_{n=1}^{N}z^{(n)}(f(x,u))r(x,u)\mathrm{d}u = \sum_{n=1}^{N} \int z^{(n)}(f(x,u))r(x,u)\mathrm{d}u = \sum_{n=1}^{N} z^{(n)}(x),
\end{align}
where we used \eqref{eq:sup_Z_rz} and \eqref{eq:sup_z_equal_Z} in the last equality. This implies that the additive structure holds at all time steps by induction.

\textbf{Computing optimal control and weights:} In order to use the additive structure in the optimal control distribution, we first write $\pi^{\star}$ in terms of $r$ and $z$ as
\begin{align*}
\pi^{\ast}(u|x) = 
\frac{1}{Z}\exp\big[ -\frac{1}{\alpha}
\big( l(x,u) + V(x')\big)
\big] = [z(x)]^{-1}r(x,u) z(f(x,u)).
\end{align*}
For the $n$ th trajectory, it has the form of 
\begin{align*}
\pi^{\ast(n)}(u|x) = [z(x)]^{-1}r(x,u) z^{(n)}(f(x,u)).
\end{align*}
Using the additive structure of $z$, and from the two equations above, $\pi^{\ast}$ is written as a weighted sum of $\pi^{\ast (n)}$ s as
\begin{align*}
\pi^{\ast}(u|x) &= [z(x)]^{-1}\Big[ \sum_{n=1}^{N}z^{(n)}(f(x,u)) \Big] r(x,u)\\
&= \sum_{n=1}^{N} \frac{z^{(n)}(x)}{z(x)}[z^{(n)}(x)]^{-1}z^{(n)}(f(x,u))r(x,u)\\
&= \sum_{n=1}^{N} \omega^{(n)}(x) \pi^{{\ast}(n)}(u|x),
\end{align*}
with the weight
\begin{align}\label{eq:sup_weight}
\omega^{(n)} &= [z(x)]^{-1}z^{(n)}(x), \quad \text{with}\ 
\sum_{n=1}^{N}\omega^{(n)} =1.
\end{align}

Each policy for the $n$ th trajectory is given by Gaussian distribution as we saw in the unimodal case:
\begin{align*}
\pi^{(n)}(\delta u^{(n)}|\delta x^{(n)}) = \mathcal{N}(\delta u^{(n)}; \delta u^{\ast (n)}, \alpha (Q_{uu}^{n})^{-1}).
\end{align*}
Since both $z$ and $\pi$ are a weighted sum of $z^{(n)}$ and $\pi^{(n)}$, the backward pass of MEDDP is achieved by performing the backward pass in $N$ different nominal trajectories and composing them with weights. To compute them, we first represent $z$ with the value function. From \eqref{eq:sup_V_islog_Z} for each trajectories, we have  
\begin{equation}\label{eq:sup_z^n_with_V^n}
z^{(n)}(x) = \exp \Big[-\frac{1}{\alpha}V^{(n)}(x)\Big].
\end{equation}
Plugging this $z^{(n)}$ in the right-hand side of  \eqref{eq:sup_additive_z}, and $Z$ in \eqref{eq:sup_V_islog_Z} in the left-hand side of the same equation, we have a similar LogSumExp structure in value function: 
\begin{align*}
 \exp\Big[ -\frac{V(x)}{\alpha}\Big] = \sum_{n=1}^{N} \exp \big[-\frac{V^{n}(x)}{\alpha}\big]
\Leftrightarrow\ \ V(x) = -\alpha \ln \sum_{n=1}^{N} \exp \Big[
-\frac{1}{\alpha}V^{(n)}(x)
\Big].
\end{align*}
Plugging the equation above in \eqref{eq:sup_V_islog_Z} gives
\begin{equation*}
z(x) = \exp \Big[ -\frac{1}{\alpha}V(x)\Big] = \sum_{n=1}^{N} \exp \Big[
-\frac{1}{\alpha}V^{(n)}(x)\Big].
\end{equation*}
From the above equation, \eqref{eq:sup_z^n_with_V^n}, and \eqref{eq:sup_weight}, we get
\begin{align*}
\omega_{t}^{(n)} &= z_{t}(x_{t})^{-1}z_{t}^{(n)}(x_{t})=\frac{\exp \big[-\frac{1}{\alpha}
(V_{t}^{(n)}(x_{t}))
\big]}
{
\sum_{n=1}^{N}\exp \big[
-\frac{1}{\alpha}
V_{t}^{(n)}(x_{t})
\big]
},
\end{align*}
where we explicitly write $V$ as a function of state and time. The value function $V_{t}(x_{t})$ can be approximated by the following quadratic approximation.
\begin{align*}
V_{t}^{(n)}(x) &\approx V_{t}^{(n)}(\bar{x}^{(n)})
+ V_{H,t}^{(n)}(\bar{x}^{(n)})
+V^{(n)}_{x,t}(\bar{x}^{(n)})^{\tr} \delta x^{(n)}
+\frac{1}{2}\delta x^{(n)\tr} V_{xx,t}^{(n)}(\bar{x}^{(n)})\delta x^{(n)},
\end{align*}
that uses derivatives from the backward pass and $\delta x$ from the forward pass. Another option to compute $V_{t}(x)$ is to compute it as the cost-to-go of a trajectory at time $t$ in the forward pass. We used the second approach in our implementation for its simplicity. 
\subsection{Algorithm}
\begin{figure}[!h]
\centering
\includegraphics[trim={2cm 9cm 5.2cm 2.5cm},clip, width=\textwidth]{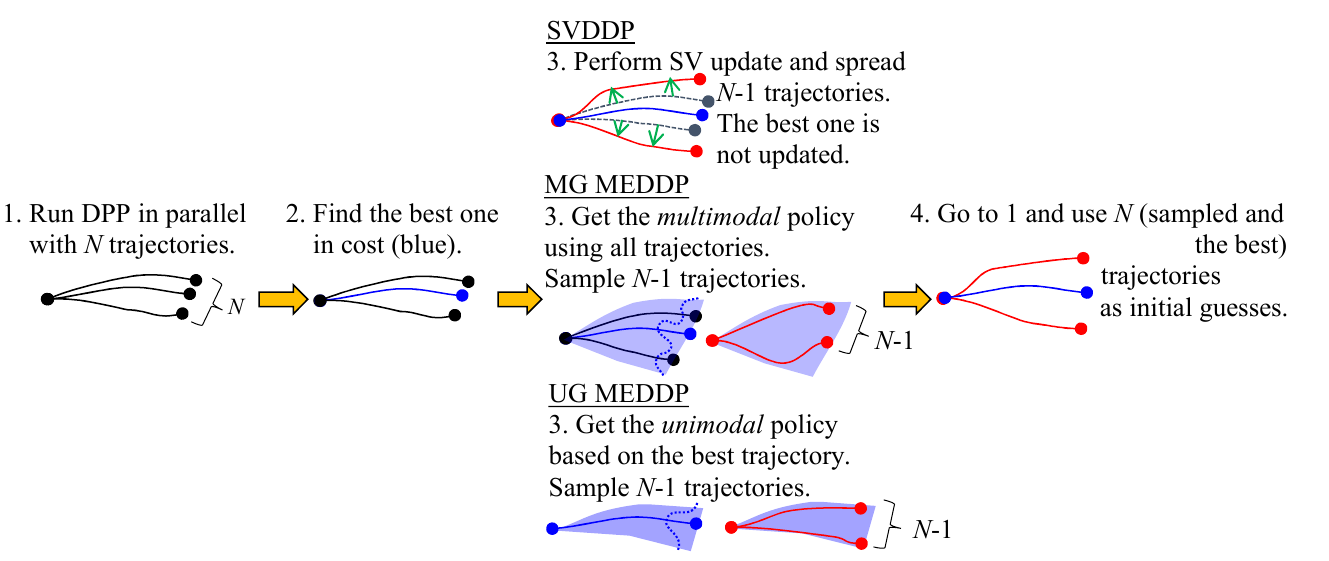}
\caption{Schematic of SVDDP, UG-MEDDP, and MG-MEDDP.}
\label{fig:sup_alg_shcematic}
\end{figure}

We provide the schematic of the MEDDP with Unimodal Gaussian (UG), Multimodal Gaussian (MG), and our SVDDP in Fig. \ref{fig:sup_alg_shcematic}. As mentioned in Section \ref{sec:SVDDP}, we preserve the best trajectory and sample the rest (in the figure, $N-1$ trajectories) from the stochastic policies in these algorithms. In MG-MEDDP, if a mode has a significantly better cost and thus better value function than others, the algorithms keep sampling from the same mode, which we call mode collapsing. When this happens, the algorithm cannot make full use of its multimodal policy. We use a heuristic to prevent this. We set a threshold for the weights or categorical probability mass and renormalize them so that the smallest weight does not fall below the threshold.


\renewcommand{\theequation}{\thesection.\arabic{equation}}
\renewcommand{\thefigure}{\thesection.\arabic{figure}} 
\renewcommand{\thetable}{\thesection.\arabic{table}}
\setcounter{figure}{0} 
\setcounter{table}{0} 
\setcounter{equation}{0}
\section{Appendix for Stein Variational DDP}

\subsection{Algorithm}\label{sec:app_svddp_code}
We provide a pseudocode of SVDDP in Algorithm \ref{alg:SVDDP} and a line search mechanism that ensures a cost in a real number after exploration in Algorithm \ref{alg:SVDDP_LS}. In these algorithms, the modes are represented by superscripts; for example, $u^{(n)}$ stands for $n$ th mode. The subscripts are for time instances. We represent a whole sequence without subscripts. For example, $u= [ u_{0}, \cdots, u_{N-1}]$. Algorithm \ref{alg:SVDDP_LS} corresponds to the selection of the step size $\epsilon$ explained in Section \ref{sec:components_SVDDP}. $\epsilon_{\rm{array}}$ given as input is the sequence of candidate of $\epsilon$s ordered in descending order ending with 0. The last 0 corresponds to not updating the trajectories with SV update. $m$ in Algorithm \ref{alg:SVDDP} is for the heuristic update (sampling) frequency explained in Section \ref{sec:components_SVDDP}. We note that although the Newton direction $w$ needs to be computed for every time step, the computation can be parallelized because no recursions are in volved.

\RestyleAlgo{ruled}
\begin{algorithm}[hbt!]
\caption{Stein Variational DDP}\label{alg:SVDDP}
\KwIn{$x_{0}$: initial condition for state, $\bar{u}$: initial control sequence, $\Sigma_{0}$: initial covariance to sample control, $N$: Number of modes, $m$: sampling frequency, $I$: max. iteration}
\SetKw{And}{and}
\SetKw{Break}{break}
\KwResult{$x^{(b)}, u^{(b)}, K^{(b)}$} 

$u^{(1)} \gets \bar{u}$\\
Sample $(N-1)$ control sequences for initial nominal trajectories. \\
\For{$n \gets 2$ to $N$ in parallel}
{
    \For{$t \gets 1$ to $T-1$}{
        $u_{t}^{(n)} \sim \mathcal{N}(\bar{u}_{t}, \Sigma_{0,t})$ 
    }
    $J^{(n)} \gets $ Compute initial cost ($x^{(n)}, u^{(n)}$)\\
    $b\gets \argmin_{n} \{J^{(n)}\}, \ n \in \{1\cdots N \}$ \  Find the best trajectory.
}
Start Optimization.\\
\For{$i \gets 0$ to $I$}
{   
    \If{$i \% m = 0$}{
        Compute coefficients for Newton's direction. See \eqref{eq:SVDDP_grad}, \eqref{eq:Hss_SVDDP}.\\
        \For{$n \gets 1$ to $N$ in parallel}{
            \For{$t \gets 0$ to $T-1$ in parallel}{
                Solve  $H^{n,n} \beta_{t}^{(n)} = -\nabla \hat{J}_{q}[\bm{0}]( u_{t}^{(n)})$
            }
        } 
    Compose Newton's direction with Galerkin approximation. See \eqref{eq:SV_ddp_update}.\\
    \For{$n \gets 1$ to $N$ in parallel}{
        \For{$t \gets 0$ to $T-1$ in parallel}{
            $w^{(n)} = \sum_{n'=1}^{N} \beta_{t}^{(n')} k(u_{t}^{(n)},u_{t}^{(n')})$}
        }
    Apply SV update except for the best one.\\
    \For{$n \gets 1$ to $N$ but not $b$ in parallel}{    
    $x^{(n)}, u^{(n)} \gets$ Line Search in Algorithm \ref{alg:SVDDP_LS} ($x^{(n)}$, $u^{(n)}$, $w^{(n)}$, $\epsilon_{\rm{array}}$).\\
    }
    }
    \For{$n \gets 1$ to $N$ in parallel}{
        $x^{(n)}, u^{(n)}, K^{(n)}, J^{(n)} \gets$ Solve one iteration of DDP ($x^{(n)}, u^{(n)}$)
        }
    $b\gets \argmin_{n} \{J^{(n)}\}, \ n \in \{1\cdots N \}$ \ Find the best trajectory.
    }
\end{algorithm}

\begin{algorithm}[hbt!]
\caption{Line Search to ensure Real and Finite Cost}\label{alg:SVDDP_LS}
\KwIn{$\bar{x}, \bar{u} , K$: sequence of state, control, and feedback gain from DDP, $w$: Newton Direction from SV, $\epsilon_{\rm{array}}$: sequence of $\epsilon$s in descending order ending with 0. }
\SetKw{And}{and}
\SetKw{Break}{break}
\SetKw{Continue}{continue}
\KwResult{$x, u$: a pair of updated trajectories.}

\For{$\epsilon$ in $\epsilon_{\rm{array}}$}
{   $x_{0} \gets \bar{x}_{0}, \ \delta{x} \gets  0$\\
    Modify control, propagate dynamics, and compute cost.\\
    \For{$t \gets 0$ to $T-1$}{
        $u_{t} \gets \bar{u}_{t} + \epsilon w_{t} + K \delta {x}$\\
        $x_{t+1} \gets f(x_{t}, u_{t})$\\
        $\delta x \gets x_{t+1} - \bar{x}_{t+1}$  
    }
    $J \gets$ Compute cost ($x,u$). (See \eqref{eq:ocp_cost_app})\\
    \uIf{$J$ is \upshape{\texttt{NaN}}}{
        \Continue \;
    }
    \Else{
    \Break\;
    }
}
\end{algorithm}

\subsection{Outer Product of Kernel}\label{sec:app_outerproduct_kernel}
This section provides an interpretation of the outer product of the gradient presented in \eqref{eq:Hss_SVDDP}.
Let us consider the case where we use Gaussian kernel $k(\theta, \theta_{n})$ and a metric $
J_{\mathrm{kernel}} = \frac{1}{2}[k(\theta, \theta_{n})]^{2}$. Taking the gradient and the Hessian gives the following terms.
\begin{align*}
\nabla_{\theta} J_{\rm{kernel}} &= k(\theta, \theta_{n}) k_{\theta}(\theta, \theta_{n}),\\
\nabla_{\theta}^{2} J_{\rm{kernel}} &= k_{\theta}(\theta, \theta_{n})k_{\theta}(\theta, \theta_{n})^{\tr} + k_{\theta} k_{\theta \theta}(\theta, \theta_{n}) \approx  k_{\theta}(\theta, \theta_{n})k_{\theta}(\theta, \theta_{n})^{\tr},  
\end{align*}
where the approximation in the Hessian is known as the Gauss-Newton approximation of the Hessian that approximates Hessian using gradient. Since the kernel is positive, the minimization of $J_{\rm{kernel}}$ is equivalent to that of
the kernel. Therefore, as claimed in the main paper, $k_{\theta}(\theta, \theta_{n})k_{\theta}(\theta, \theta_{n})^{\tr}$ captures information on in which direction and how much the kernel changes with respect to $\theta$.

\subsection{Kernel and Repulsive Force}\label{sec:app_kernel_repulsive_force}
This section provides the details of the kernel that we used in the main paper. We use the Gaussian kernel with the following form:
\begin{align*}
k(\theta, \theta_{n}) = \exp(-\norm{\theta-\theta_{n}}^{2})/L, \quad L = \text{med}\{ \norm{\theta-\theta_{n}}^{2} \}/\log(d), \quad \text{for} \ n=1, \cdots, N, \  
\end{align*}
with $\theta \in \mathbb{R}^{d}$. Here, the length parameter $L$ is chosen with the median heuristic of the squared distance to all data points. In MPC with SVGD, the optimization parameter was mean control trajectories, leading to high dimensional $\theta$ which has (dimension of control) $\times$ (prediction horizon of MPC) entries \citep{Wang2021Tsallis, lambdert2021SVMPC}. As a result, SVGD becomes less effective with a small repulsive force. They used a sum of local kernels to alleviate this issue. In contrast, in our work, the parameter $\theta$ only has a dimension of control due to the nature of DDP, i.e., splitting a problem associated with a trajectory into small subproblems at every time step. Therefore, our algorithm works well with a normal Gaussian kernel in contrast to the existing SVGD MPC methods.

\begin{figure}[ht]
\centering
\includegraphics[width=\textwidth]{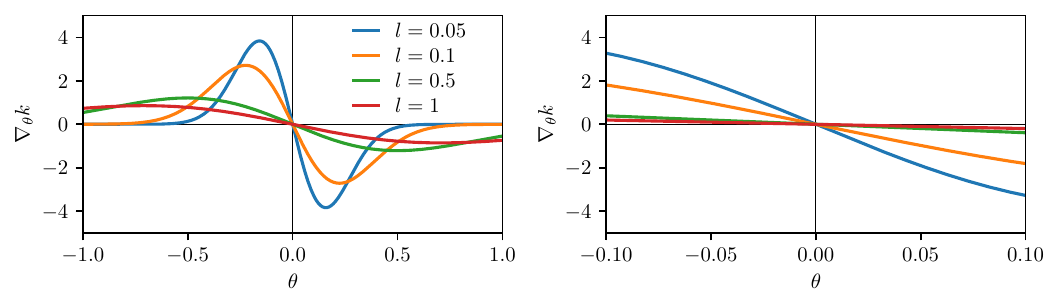}
\caption{Relationship between $\theta$ and $\nabla_{\theta}k(\theta, \theta')$ with $\theta'=0$.}
\label{fig:dkernel}
\end{figure}

To understand the behavior of the repulsive force governed by $\nabla_{\theta}k(\theta, \theta')$, 
We visualize $\nabla_{\theta}k(\theta, \theta')$ with a scalar variable $\theta$, a fixed value $\theta'(=0)$, and different values of fixed $L$ in Fig. \ref{fig:dkernel}. Let us observe the change of gradient term. We start $\theta$ that is far from $\theta'$ and make $\theta$ close to $\theta'=0$. The repulsive force increases as $\theta$ approaches $\theta'(=0)$ and shows its maximum. After that, it decreases again and becomes zero at $\theta=\theta'(=0)$. To see the effect of the median heuristic, we consider the $N$ number of $\theta'$s denoted by $\theta_{n}$, $n=1,\cdots, N$. With the median heuristic, when many $\theta_{n}$s exist that take a close value, and when $\theta$ also becomes close to that value, $L$ becomes small. This adjusted $L$ prevents the repulsive force from vanishing by scaling the squared distance part of the Gaussian kernel. We note that if $\theta$ becomes identical to $\theta'$, the repulsive force vanishes.

\renewcommand{\theequation}{\thesection.\arabic{equation}}
\renewcommand{\thefigure}{\thesection.\arabic{figure}} 
\renewcommand{\thetable}{\thesection.\arabic{table}}
\setcounter{figure}{0} 
\setcounter{table}{0} 
\setcounter{equation}{0}
\section{DDP with Barrier Term}\label{section:app_DDP_with_barrier}
In this section, we provide the derivation of constrained DDP with the penalty function introduced in section \ref{sec:handle_const_DDP}. We first introduce a relaxed version of the $\log$ barrier function. Then, we derive the constrained DDP with the barrier term. Finally, we explain the conditioning of the Hessian.

\subsection{Barrier methods and Relaxed-Barrier Function}
Consider a constrained optimization problem with the same cost in \eqref{eq:ocp_cost} under constraints: 
\begin{align}\label{eq:const_inequality}
g_{t}(x_{t}, u_{t}) &\leq 0,  \quad \text{for} \quad t=0  \cdots T-1, \quad \text{and}\quad
g_{T}(x_{T}) \leq 0.
\end{align}
with $g_{t}\in \Rb^{w_{1}}$ and $g_{T}\in \Rb^{w_{2}}$. To simplify the notation, we concatenate the state and control vector, creating a new vector as follows.
\begin{align*}
y_{t} = [x_{t}^{\tr}, u_{t}^{\tr}]^{\tr}, \quad y_{T} = x_{T}.
\end{align*}
Barrier methods solve this problem by minimizing the sequence of a new objective which is the sum of the original objective and the constraints transformed by barrier functions parameterized by a nonnegative penalty (barrier) parameter $\mu$ \cite{Frisch1955logarithmic,Fiacco1968sequential}. We take a logarithmic barrier function, and define a penalty term by $\mathcal{P}$ as
\begin{align}{\label{eq:log_barrier_penalty_function}}
    \mathcal{P}(y_{t};\mu) = 
    \begin{cases}   
     - \mu \sum_{i=1}^{w_{0}}\log(-g_{i}(y_{t})), & g(y_{t}) < 0,\\
    \infty, & g(y_{t}) \geq  0,
    \end{cases}
\end{align}
where we use $w_{0}$ to represent $w_{1}$ and $w_{2}$ to avoid excessive notation. Due to the $\log$ term, the penalty term increases as $y$ approaches the boundary of $g(y)=0$ and explodes when $g(y) > 0$. Typically, the new objective (the sum of the original and the penalty term) is minimized with a fixed $\mu$, followed by a reduction of $\mu$ to bring the minimizer of the new objective closer to that of the original objective \cite{Nocedal2006numerical}. However, fixed $\mu$ can work well in practice \cite{Grandia2019BarrierDDP}. By ensuring the finiteness of $\mathcal{P}$, the constraints are satisfied for all iterations, even before convergence.

One of the main drawbacks of the $\log$ barrier method is that it cannot accept an infeasible trajectory, which does not fit well with exploration in our setting. Because the sampling process is not explicitly aware of the constraints, the sampled trajectories for DDP can be infeasible. To solve the problem, a smooth quadratic approximation (relaxation) of the $\log$ barrier function \cite{Hauser2006relaxbarrier} can be used, which is given by  
\begin{align}\label{eq:relaxed_barrier}
\mathcal{P}_{\rm{r}}(y;\mu) = 
\begin{cases}
-\mu \sum_{i=1}^{w_{0}}\log (-g_{i}(y)),  \quad \quad  & g(y) \leq -\delta, \\
\mu \sum_{i=1}^{w_{0}}\Bigg[ \frac{1}{2}\Big[\Big(\frac{g_
{i}(y)+2\delta}{\delta}\Big)^{2} - 1\Big] - \log{\delta} \Bigg],   & g(y) > -\delta.
\end{cases}
\end{align}
\subsection{Constrained DDP}
We now solve \eqref{eq:ocp_cost} under constraints $g(y_{t})$. By adding \eqref{eq:relaxed_barrier} to the cost, we have a modified cost $\hat{J}$ as
\begin{align}\label{eq:log_barr_DDP_problem}
\hat{J}(X,U) &= \sum_{t=0}^{T-1}\Big[ l(y_{t})  {\color{red}{+ \mathcal{P}_{{\rm{r}}}(y_{t})}}{\color{black}\Big]} + \Phi(y_{N})+ {\color{red}{\mathcal{P}_{{\rm{r}}}(y_{T})}},
\end{align}
which modifies the Bellmas's optimality principle as follows.
\begin{align*}
    {V}_{t}(x_{t}) = \min_{u_{t}}\big[    
    \underbrace{l(y_{t}) 
    {\color{red}+ \mathcal{P}_{\rm{r}}(y_{t})}}_{\hat{Q}(y_{t})}+ {V}_{t+1}(x_{t+1})\big].
\end{align*}
Therefore, the problem is solved with DDP by changing $Q$ to $\hat{Q}$.
\subsection{Conditioning Hessian}
In DDP, maintaining $Q_{uu}$ PD is critical to ensure cost reduction and convergence. In this section, we examine to what extent the added penalty term affects the PD property of $\hat{Q}_{uu}$. Consider the relaxed barrier in \eqref{eq:relaxed_barrier}. When $g(y)> -\delta$, the Hessian of $\mathcal{P}_{\rm{r}}(y)$ is also PD. Therefore, we only need to ensure that $Q_{uu}$ is PD as in normal (unconstrained) DDP. However, when $g(y)\leq -\delta$, we have
\begin{align*}
    \hat{Q}_{yy} &= Q_{yy} {\color{red}{- \mu \sum_{i=1}^{w}\left[\frac{g_{i,yy}(y)}{g_{i}(y)}\right] + \mu \sum_{i=1}^{w}\Big[\frac{g_{i,y}(y)g_{i,y}(y)^{\tr}}{g_{i}(y)^{2}}\Big]}},
\end{align*}
where the third term is PSD by construction, but the second term might not be. 
The second term may slow down optimization and even make it fail.
Therefore, we exclude the second term for the better conditioning of $\hat{Q}_{uu}$. This approximation corresponds to the Gauss-Newton approximation in the optimization technique, where the Hessian is approximated based on the gradient \cite{Nocedal2006numerical}.

\renewcommand{\theequation}{\thesection.\arabic{equation}}
\renewcommand{\thefigure}{\thesection.\arabic{figure}} 
\renewcommand{\thetable}{\thesection.\arabic{table}}
\setcounter{figure}{0} 
\setcounter{table}{0} 
\setcounter{equation}{0}

\section{Results on Trajectory Optimization}\label{sed:appendix_TO}
This section provides the results of the Trajectory Optimization (TO) that is referred to in the main paper. 
The result includes the figures for the output and the mean and standard deviation of the cost with different temperatures $\alpha$ for ME and SVDDPs and $\epsilon$ for SVDDP. We also give a belief explanation of the experimental setup, such as cost structures.

\subsection{Overview of the Experiment}\label{sec:app_TO_overview}
\begin{figure}[t]
\begin{subfigure}[b]{0.25\linewidth}
\includegraphics[trim={0.3cm 0.3cm 0.3cm 1.0cm},clip,width=\textwidth]{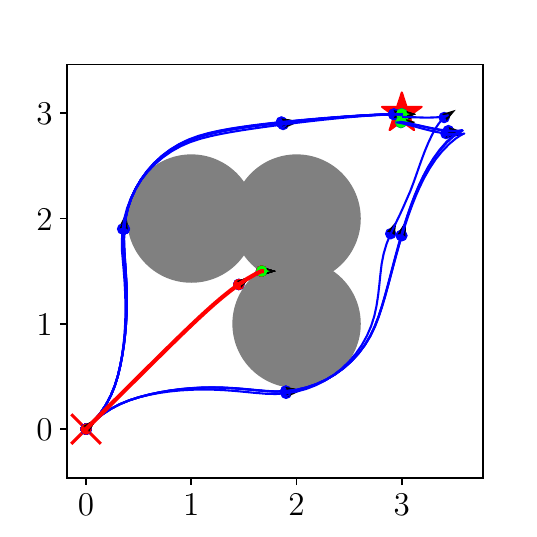}
\caption{2D Car}
\end{subfigure}    
\begin{subfigure}[b]{0.25\linewidth}
\centering
\includegraphics[trim={1.5cm 1.7cm 2cm 2.5cm},clip,width=\textwidth]{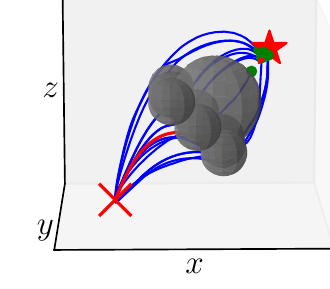}
\caption{Quadrotor}
\end{subfigure}
\begin{subfigure}[b]{0.48\linewidth}
\centering
\includegraphics[trim={0.3cm 0.3cm 0.3cm 1.05cm},clip,width=\textwidth]{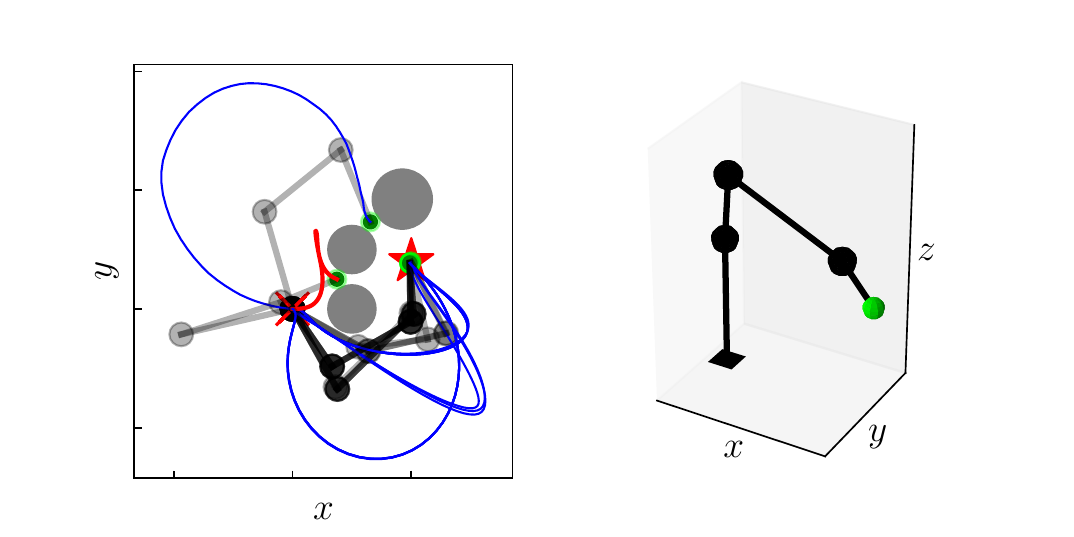}
\caption{Arm}
\end{subfigure}
\caption{TO Tasks for the three systems.{\color{red}$\times$} and {\color{red}$\bigstar$} indicate the initial and target positions. The positions at the terminal time step in the 2D car and quadrotor and that of the arm's end effector are plotted in green. A trajectory that is obtained by normal (no exploration) DDP is drawn in red. The arm is projected on to $x$-$y$ plane.}
\label{fig:TO_overview}
\end{figure}

Fig. \ref{fig:TO_overview} provides an overview of the TO experiments with three systems, i.e., 2D car, quadrotor, and robotic arm. In the experiments, we use the following quadratic cost with the target state $x_{\rm{g}}$.
\begin{equation}
    J(X,U) = \sum_{t=0}^{T-1} \Big[ \frac{1}{2} u_{t}^{\tr}Ru_{t} +  \frac{1}{2}(x_{t}-x_{\rm{g}})^{\tr} P (x_{t}-x_{\rm{g}}) \Big] + \frac{1}{2}(x_{T}-x_{\rm{g}})^{\tr}Q(x_{T}-x_{\rm{g}}),
\end{equation}
where $P$, $Q$, and $R$ are diagonal weight matrices for running state cost, terminal (sate) cost, and control cost. Based on this cost, constraint terms are added to satisfy constraint $g(x_{t},u_{t})\leq0$. 
penalty term for constraints given in \eqref{eq:log_barr_DDP_problem} is added:
\begin{align*}
\hat{J}_{\rm{DDP}} = J(X,U) + \Big[\sum_{t=0}^{T-1}\mathcal{P}(g(x_{t},u_{t}))\Big] + \mathcal{P}(g(x_{T})). 
\end{align*}

\subsection{Mean and Standard Deviation of Cost}
The mean and standard deviation of the TO experiment are provided in Table \ref{tab:TO_alpha_mode}. The best mean cost among the same number of modes and $\alpha$ is in bold.

\begin{table}[h!]
\footnotesize
\caption{Results of Trajectory Optimization.}
\label{tab:TO_alpha_mode}
\begin{subtable}[t]{\textwidth}
\caption{2D Car TO Result.}\label{tab:car_alpha_mode}
\resizebox{\columnwidth}{!}{
\begin{tabular}{lllllllllllllll}
\hline
\multicolumn{1}{|l|}{\multirow{2}{*}{$\alpha$}} &
  \multicolumn{1}{l|}{\multirow{2}{*}{num}} &
  \multicolumn{2}{l|}{Uni.} &
  \multicolumn{2}{l|}{Multi.} &
  \multicolumn{2}{l|}{SV ($\epsilon=1$)} &
  \multicolumn{2}{l|}{SV ($\epsilon=2$)} &
  \multicolumn{2}{l|}{SV ($\epsilon=3$)} &
  \multicolumn{2}{l|}{SV ($\epsilon=4$)}
  \\ \cline{3-14} 
\multicolumn{1}{|l|}{} &
  \multicolumn{1}{l|}{} &
  \multicolumn{1}{l|}{Mean} &
  \multicolumn{1}{l|}{Std.} &
  \multicolumn{1}{l|}{Mean} &
  \multicolumn{1}{l|}{Std.} &
  \multicolumn{1}{l|}{Mean} &
  \multicolumn{1}{l|}{Std.} &
  \multicolumn{1}{l|}{Mean} &
  \multicolumn{1}{l|}{Std.} &
  \multicolumn{1}{l|}{Mean} &
  \multicolumn{1}{l|}{Std.} &
  \multicolumn{1}{l|}{Mean} &
  \multicolumn{1}{l|}{Std.}
  \\ \hline
\multirow{4}{*}{10} 
& 8  & 306 & 115  & \textbf{250}   & 82    & 405 & 85.3  & 417 & 78.2 & 311 & 113 & 255 & 82.3\\
& 16 & 290  & 111 & 280   & 102   & 354 & 80.4 & 270 & 93.0 & \textbf{211} & 0.498 & 225 & 37.8\\
& 32 & 224 & 48.6 & 212   & 0.389 & 212 & 113 & \textbf{211} & \num{7.39e-2} & \textbf{211 }& 0.219 & \textbf{211 }& 1.17\\
& 64 & \textbf{211} & 0.626 & 212   & 0.572 & 417 & 3.99 & \textbf{211} & \num{218e-2} & \textbf{211} & 0.0279  & \textbf{211} & 0.0377 \\ \hline
\multirow{4}{*}{20}
& 8  & 275 & 104    & \textbf{274 }  & 104   & 400 & 94.3 &309 & 113 &263 & 83.4 &300 & 107  \\
& 16 & 274 & 104    & \textbf{211}   & 0.710 & 447 & 9.95 & 217 & 20.6 & \textbf{211} & 0.430 & 215 & 13.3   \\
& 32 & 216 & 0.711  & 212   & 0.645 & 301 & 107 & \textbf{211} & 0.234 & \textbf{211 }& 0.974 & \textbf{211} & 1.21\\
& 64 & \textbf{211} & 0.741 & \textbf{211}   & 0.712 & \textbf{211} & \num{6.13e-2} & \textbf{211} & 0.379 & \textbf{211} & 0.682 & 212 & 1.99  \\ \hline
\multirow{4}{*}{30}
& 8  & 369  & 111 & \textbf{243}   & 80.1  & 405 & 85.4 & 263 & 93.8 & 249 & 81.4 & 265 & 93.8 \\
& 16 & 296 & 109  & 227   & 58.8  & 411 & 80.2 & \textbf{211 }& 1.22 & \textbf{211} & 0.613 & 256 & 90.2 \\  
& 32 & \textbf{211} & 0.707 & \textbf{211} & 0.628 & 290 & 111 & \textbf{211} & 0.154 & 212 & 1.01 & 217 & 19.5 \\
& 64 & \textbf{211} & 0.694 & 212 & 0.647 & \textbf{211} & 0.0185 & \textbf{211} & 2.14 & \textbf{211} & 0.592 & 212 & 0.727 \\ \hline
\multirow{4}{*}{40}
& 8  & 389  & 97.3           & 269 & 70.0  & 400 & 94.0 & 242 & 67.1 & 273 & 94.7 & \textbf{234} & 44.8  \\
& 16 & 329 & 110 & 227   & 58.8  & 396 & 87.6 & \textbf{211} & 0.9 & 241 & 77.1 & \textbf{211} & 0.896\\
& 32 & 247  & 80   & 212 & 0.751 & \textbf{211} & 0.907 & 212 & 1.54 & \textbf{211} & 0.28 & \textbf{211} & 0.482\\            
& 64 & \textbf{211} & 0.649 & \textbf{211}   & 0.714 & \textbf{211} & 0.0924 & \textbf{211} & 0.626 & 212 & 4.55 & \textbf{211} & 0.566\\\hline
\end{tabular}}
\end{subtable}

\begin{subtable}[t]{\textwidth}
\caption{Quadrotor TO Result.}
\resizebox{\columnwidth}{!}{
\begin{tabular}{llllllllllllll}
\hline
\multicolumn{1}{|c|}{\multirow{2}{*}{$\alpha$}}& 
\multicolumn{1}{c|}{\multirow{2}{*}{\# of modes}} & 
\multicolumn{2}{l|}{Uni.} &
\multicolumn{2}{l|}{Multi.} &
\multicolumn{2}{l|}{SV ($\epsilon=1$)} &
\multicolumn{2}{l|}{SV($\epsilon=2$)} &
\multicolumn{2}{l|}{SV ($\epsilon=3$)} &
\multicolumn{2}{l|}{SV($\epsilon=4$)}\\ \cline{3-14} 
\multicolumn{1}{|c|}{} &
\multicolumn{1}{c|}{} &
\multicolumn{1}{l|}{Mean} &\multicolumn{1}{l|}{Std.} &
\multicolumn{1}{l|}{Mean} & \multicolumn{1}{l|}{Std.} &
\multicolumn{1}{l|}{Mean} & \multicolumn{1}{l|}{Std.} &
\multicolumn{1}{l|}{Mean} & \multicolumn{1}{l|}{Std.} & 
\multicolumn{1}{l|}{Mean} & \multicolumn{1}{l|}{Std.}&
\multicolumn{1}{l|}{Mean} & \multicolumn{1}{l|}{Std.}
\\ \hline
\multirow{4}{*}{0.01} 
& 8
& 307 & 79.4 & \textbf{146} & 78.4 & 291 & 90.2 & 260 & 106 & 169 & 102 & 152 & 92.2 \\
& 16
& 184 & 110 & \textbf{{107}} & 6.26 & 290 & 91.0 & 181 & 109 & 167 & 102 & 135 & 78.5 \\
&32
& 153 & 92.8 & 109 & 5.51 & 259 & 106 & 135 & 78.5 & 104 & 5.55 & \textbf{103} & 4.61 \\
&64
& 155 & 92.3 & 108 & 4.63 & 183 & 107 & 118 & 58.0 & \textbf{{101}} & 0.0198 & \textbf{{101}} & 0.0165 \\ \hline
\multirow{4}{*}{0.1} & 8 & 159 & 91.0 & 111 & 7.14 & 245 & 112 & 108 & 6.35 & 105 & 3.98 & \textbf{106} & 5.76 \\
 & 16 & 111 & 4.80 & 107 & 6.77 & 184 & 107 & \textbf{104 }& 4.19 & \textbf{104} & 4.3 & \textbf{104} & 3.96 \\
 & 32 & 106 & 6.18 & 106 & 6.49 & 215 & 112 & 106 & 4.68 & \textbf{103} & 3.42 & \textbf{103} & 3.4 \\
 & 64 & 108 & 5.99 & 104 & 5.63 & 106 & 4.39 & \textbf{103} & 3.41 & \textbf{103} & 3.39 & \textbf{101} & 0.0168 \\ \hline
\multirow{4}{*}{1} & 8 & 186 & 107 & 109 & 6.22 & 124 & 56.7 & 107 & 5.05 & \textbf{105} & 5.43 & 107 & 4.58 \\
 & 16 & 162 & 94.8 & 111 & 7.90 & 104 & 4.01 & 104 & 3.74 & 103 & 3.27 & \textbf{102} & 3.53 \\
 & 32 & 111 & 4.76 & 109 & 7.20 & 105 & 4.26 & 104 & 4.08 & \textbf{101} & 0.537 & 102& 2.93 \\
 & 64 & 110 & 4.26 & 111 & 5.32 & 104 & 3.93 & \textbf{101} & 0.182 & \textbf{101} & 0.0327 & \textbf{101} & 2.14  \\ \hline
\multirow{4}{*}{5}
& 8 & 270 & 114 & 210 & 116 & 108 & 4.69 & \textbf{106} & 4.98 & 109 & 5.74 & 124 & 57.3 \\
 & 16 
 & 186 & 107 & 112 & 5.40 &106 & 4.16 & \textbf{104} & 3.37 & 106 & 3.56 & 107 & 3.64  \\
 & 32 
 & 158 & 90.7 & 110 & 6.62 & 103 & 3.3 & \textbf{102} & 2.08 & 103 & 2.58 & \textbf{102} & 1.1 \\
 & 64 
 & 125 & 57.0 & 111 & 4.63 & 103 & 3.36 & 102 & 2.09 & \textbf{101} & 0.78 & 102 & 1.23\\ \hline
\multirow{4}{*}{10} 
& 8 
& 240 & 121 & 177 & 110 & \textbf{109} & 5.84 & \textbf{109} & 2.15 & 110 & 5.01 & 111 & 4.39  \\
& 16 
& 216 & 115 & 146 & 78.8 & \textbf{106 }& 3.78 & 107 & 4.12 & 109 & 4.95 & 109 & 4.74 \\
&32
& 125 & 57.1 & 110 & 2.09 & \textbf{103 }& 2.36 & \textbf{103} & 2.65 & 107 & 3.14 & 106 & 3.63\\
&64
& 109 & 1.84 & 109 & 2.53 & 104 & 3.83 & 103 & 2.99 & \textbf{102} & 2.32 & \textbf{102} & 1.42 \\
\hline
\end{tabular}}
\end{subtable}

\begin{subtable}[t]{\textwidth}
\caption{Arm TO Result.}
\resizebox{\columnwidth}{!}{
\begin{tabular}{llllllllllllll}
\hline
\multicolumn{1}{|c|}{\multirow{2}{*}{$\alpha$}}& 
\multicolumn{1}{c|}{\multirow{2}{*}{\# of modes}} & 
\multicolumn{2}{l|}{Uni.} &
\multicolumn{2}{l|}{Multi.} &
\multicolumn{2}{l|}{SV ($\epsilon=1$)} &
\multicolumn{2}{l|}{SV($\epsilon=2$)} &
\multicolumn{2}{l|}{SV ($\epsilon=3$)} &
\multicolumn{2}{l|}{SV($\epsilon=4$)}\\ \cline{3-14} 
\multicolumn{1}{|c|}{} &
\multicolumn{1}{c|}{} &
\multicolumn{1}{l|}{Mean} &\multicolumn{1}{l|}{Std.} &
\multicolumn{1}{l|}{Mean} & \multicolumn{1}{l|}{Std.} &
\multicolumn{1}{l|}{Mean} & \multicolumn{1}{l|}{Std.} &
\multicolumn{1}{l|}{Mean} & \multicolumn{1}{l|}{Std.} & 
\multicolumn{1}{l|}{Mean} & \multicolumn{1}{l|}{Std.}&
\multicolumn{1}{l|}{Mean} & \multicolumn{1}{l|}{Std.}
\\ \hline
\multirow{4}{*}{10}
&8
&148 & 42.1 & 133 & 31.6 & 143 & 39.8 & 149 & 37.9 & \textbf{129} & 18.5 & 134 & 37.2 \\
&16
&119 & 18.9 & 120 & 19.0 & 124 & 17.6 & \textbf{112} & 10.4 & 123 & 24.6 & 128 & 21.6 \\ 
&32
&124 & 27.0 & 113 & 5.98 & \textbf{111} & 5.65 & 113 & 5.36 & 109 & 6.68 & \textbf{111} & 7.66 \\
&64
&116 & 16.1 & 113 & 7.18 & 107 & 4.88 & 109 & 3.61 & \textbf{111} & 1.56 & 113 & 6.39 
\\ \hline
\multirow{4}{*}{20}
&8
&144 & 40.4 & \textbf{119} & 25.4 & 135 & 36.9 & 140 & 40.6 & 138 & 37.1 & 130 & 31.3  \\
&16
&123 & 20.6 & 123 & 29.6 & \textbf{111} & 6.6 & 119 & 17.0 & 113 & 6.25 & 119 & 8.43  \\
&32
&117 & 12.2 & 115 & 8.94 & \textbf{111} & 18.5 & 109 & 6.25 & 115 & 5.96 & 116 & 7.26  \\
&64
&112 & 6.77 & 113 & 5.95 & \textbf{108} & 4.72 & 114 & 4.87 & 113 & 5.09 & 114 & 5.38   \\ \hline
\multirow{4}{*}{30}
&8
&130 & 35.1 & 133 & 35.8 & 133 & 36.6 & 157 & 45.6 & 134 & 34.9 & \textbf{128} & 29.0  \\
&16
&112 & 6.3 & 115 & 10.1 & 129 & 36.4 & 117 & 7.15 & 116 & 8.71 & \textbf{113} & 7.74  \\
&32
&113 & 5.69 & 113 & 4.38 & \textbf{108} & 5.98 & 114 & 6.55 & 117 & 7.58 & 117 & 7.52  \\
&64
&115 & 6.72 & 114 & 4.76 & \textbf{111} & 0.905 & 112 & 6.14 & 114 & 4.06 & 113 & 9.16  \\ \hline
\multirow{4}{*}{40} 
&8
&130 & 35.2 & 127 & 30.2 & 129 & 30.8 & 131 & 29.2 & 129 & 26.2 & \textbf{125} & 26.3  \\
&16
&133 & 36.6 & 112 & 6.57 & 124 & 23.5 & \textbf{111} & 4.95 & 123 & 16.3 & 119 & 14.7  \\
&32
&114 & 5.16 & \textbf{113} & 5.9 & \textbf{113} & 15.5 & 116 & 5.14 & 115 & 6.71 & 115 & 6.48  \\
&64
&114 & 6.31 & 113 & 4.38 & \textbf{111} & 5.57 & 116 & 5.46 & 113 & 6.41 & 113 & 3.71  \\ \hline
\end{tabular}}
\end{subtable}
\end{table}

\begin{figure}[h]
\begin{subfigure}[b]{0.3\linewidth}
\includegraphics[trim={0.5cm 0.3cm 0.4cm 0.3cm},clip,width=\textwidth]{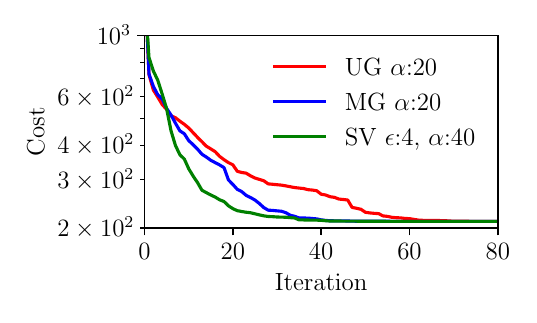}
\caption{2D Car}
\label{}
\end{subfigure}    
\begin{subfigure}[b]{0.32\linewidth}
\centering
\includegraphics[trim={0.3cm 0.3cm 0.3cm 0.3cm},clip,width=\textwidth]{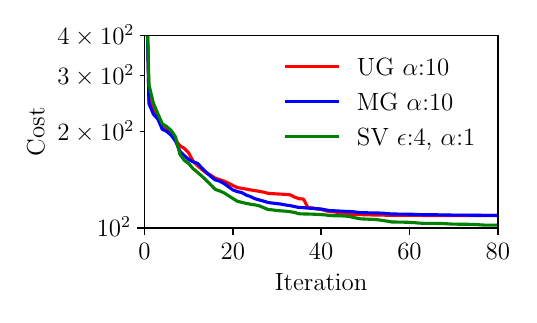}
\caption{Quadrotor}
\label{fig:}
\end{subfigure}
\begin{subfigure}[b]{0.32\linewidth}
\centering
\includegraphics[trim={0.3cm 0.3cm 0.3cm 0.3cm},clip,width=\textwidth]{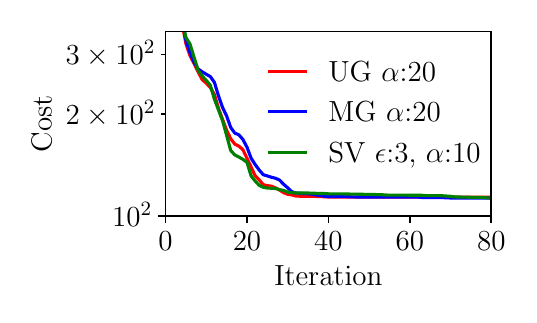}
\caption{Arm}
\label{fig:car_uni_me_ddp_mpc}
\end{subfigure}
\caption{Convergence of UG-MEDDP, MG-ME, and SVDDPs.}
\label{fig:TO_convergence}
\end{figure}

\subsection{Convergence Speed}
Fig. \ref{fig:TO_convergence} provides graphs of the evolution of the cost over iterations. The graphs are drawn with the best trajectories in terms of the convergence speed. The parameters, the temperature $\alpha$ for MEDDPs, and $\alpha$ and the step size $\epsilon$ for SV DDP, that give these results are provided in the graphs. All results here are with 64 modes.

\section{Results on Model Predictive Control}\label{sec:app_MPC_results}
This section further details the Model Predictive Control example given in the main paper. 

\subsection{2D Car} \label{sec:app_MPC_2D_car}
In Section \ref{MPC_experiments}, we evaluate the performance of SVDDP, MG-MEDDP, UG-MEDDP, and DDP in MPC mode by the 2D car dynamics. These algorithms are equipped with barrier terms to handle constraints. The algorithms, except for DDP were tested ten times to capture their stochasticity in sampling. The cost structure used in the experiment is the same as that of TO. We set the time discretization interval $dt=0.02$ s and the prediction horizon of $T_{\rm{predict}}=60$ steps. The total duration of the experiment was $T=200$ steps. For the SV and MEDDPs number of modes (trajectories) optimized in parallel was eight. For UG-MPPI and MG-MPPI, we use 2048 samples. For MG-MPPI, the number of modes is 8. For SV-MPPI, the number of Stein particles is 8, each of which has 256 samples.


\subsection{Quadrotor} \label{sec:app_MPC_Quadrotor} 
The time discretization interval was chosen as $dt = 0.01$s, and the prediction horizon was $T_{\rm{predict}} = 50$ steps. The total time step of the experiment was $T = 350$ steps. 
Other settings are the same as that of the 2D car experiment.
For UG-MPPI and MG-MPPI, we use 8000 samples. For MG-MPPI, the number of modes is 8. For SV-MPPI, the number of Stein particles is 8, each of which has 1000 samples.

\subsection{Smoothness of Control:} \label{sec:app_smooth_ctrl}DDPs can generate smoother control trajectories than MPPIs as presented in Figure \ref{fig:mpc_control_cmp}. In this figure, the dotted lines indicate the control limit, showing that controls are inside bounds. 



\begin{figure}[h]
\begin{subfigure}{0.5\textwidth}\centering
\includegraphics[trim={0.0cm 3.2cm 0.0cm 0.12cm},clip,width=1\textwidth]{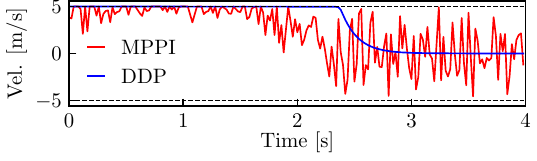}
\caption{2D car}
\end{subfigure}\hfill
\begin{subfigure}{0.5\textwidth}
\includegraphics[trim={0.0cm 0.5cm 0.0cm 2.7cm},clip,width=1\textwidth]{fig/MPC/quad/Paper/Comparison_Control_SV_stacked.pdf}
\caption{Quadrotor}
\end{subfigure}
\caption{Control trajectories of SVDDP MPC and SV-MPPI. Dotted lines indicate control limits.}
\label{fig:app_cntl_smoothness}
\end{figure}

\subsection{MPPI for Comparison}
Finally, we briefly review the cost structure of MPPIs in our comparison. The MPPIs in our experiment are implemented based on \cite{Wang2021Tsallis}.
In MPPI, a high crush cost is added to penalize state constraint violation to the original trajectory cost:
\begin{align*}
\hat{J}_{\rm{MPPI}} = J(X,U) + \Big[\sum_{t=0}^{T}c_{\rm{crush}}\mathbbm{1}_{\tilde{g}(y_{t}) > 0} \Big],  
\end{align*}
where $\mathbbm{1}_{\tilde{g}(x_{t})}$ is an indicator such that 
\begin{align*}
\mathbbm{1}_{\rm{condition}} = 
\begin{cases}
1, \quad \text{if \ condition is True,} \\
0, \quad \text{otherwise}.
\end{cases}
\end{align*}
we divide the constraints into state-dependent parts, e.g., obstacles to avoid, denoted by $\tilde{g}(x_{t})$, and control-dependent parts. The control-dependent part is satisfied by clamping the sampled control trajectories that exceed the limits. The forward propagation of the dynamics stops once the system touches the obstacles or the constraints are activated.

\section{Additional Experiments}\label{sec:add_ant_swimmer}
This section provides additional experimental results with the ant and swimmer dynamics implemented in Brax to see our algorithms' performance in underactuated and complex nonlinear systems with multibody dynamics. The experiments  are organized into  two sections. In the first section, we vary the temperature parameter $\alpha$, and let the algorithm solve reaching tasks to observe the effect of $\alpha$. In the second section, we use the ant dynamics and run similar experiments in five different obstacle fields to compare the performance of the algorithms. 

\begin{figure}[t]
\centering
\begin{subfigure}[b]{1\linewidth}
\includegraphics[width=\textwidth]{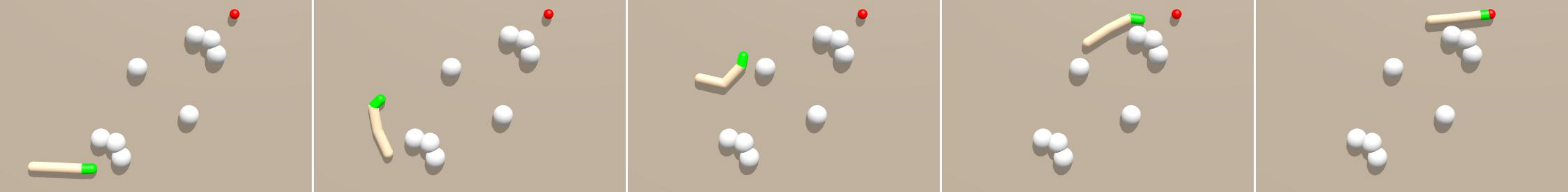}
\caption{Success}
\label{fig:swimmer_over_1}
\end{subfigure}
\begin{subfigure}[b]{1\linewidth}
\centering
\includegraphics[width=\textwidth]{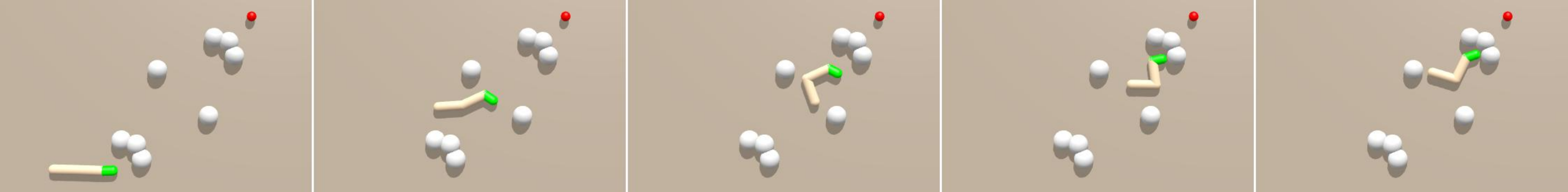}
\caption{Failure}
\label{fig:swimmer_over_2}
\end{subfigure}
\caption{Examples of a success and a failure (poor local minimum) of the swimmer example. The obstacles and targets are white and red spheres, respectively. The link with nose is shown in green. Time elapses from left to right.}
\label{fig:swimmer_overview}
\end{figure}

\begin{figure}[t]
\centering
\begin{subfigure}[b]{1\linewidth}
\includegraphics[width=\textwidth]{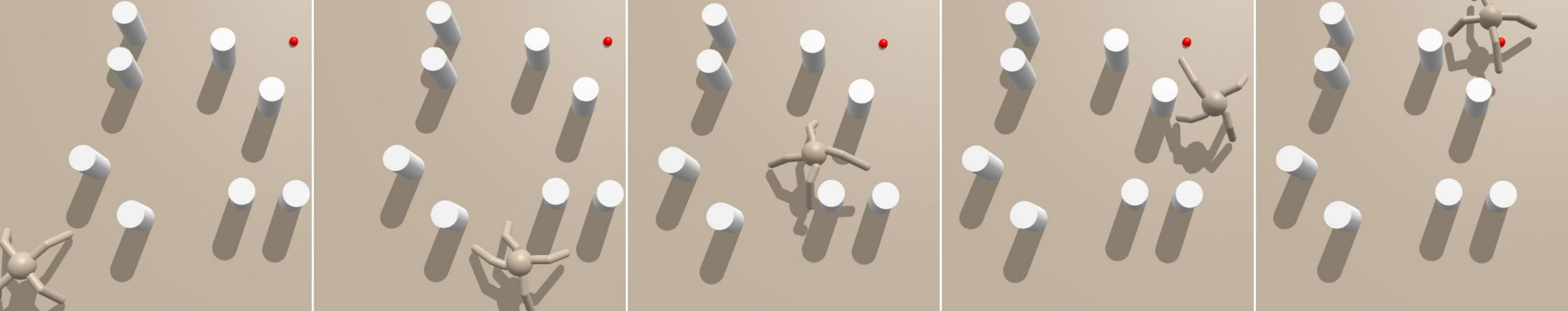}
\caption{Success}
\label{fig:ant_over_1}
\end{subfigure}
\begin{subfigure}[b]{1\linewidth}
\centering
\includegraphics[width=\textwidth]{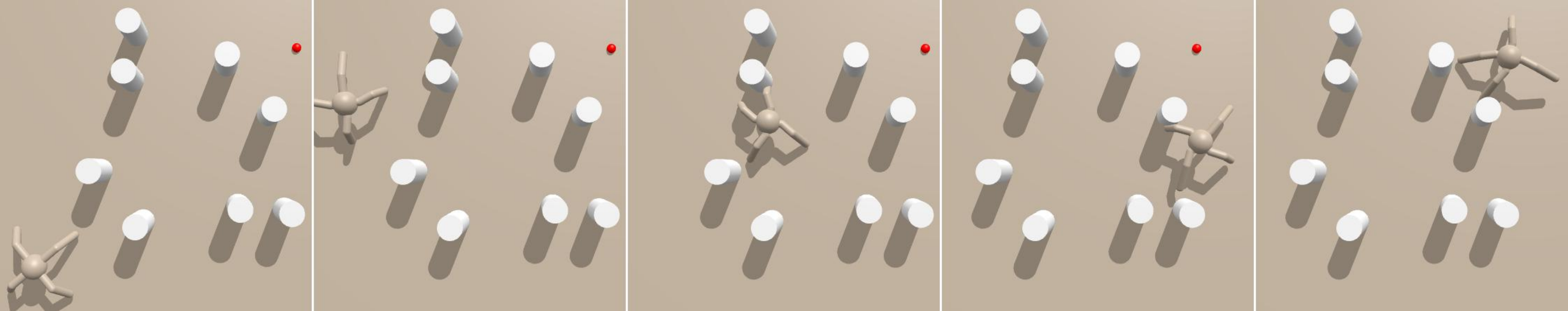}
\caption{Success}
\label{fig:ant_over_2}
\end{subfigure}
\begin{subfigure}[b]{1\linewidth}
\centering
\includegraphics[width=\textwidth]{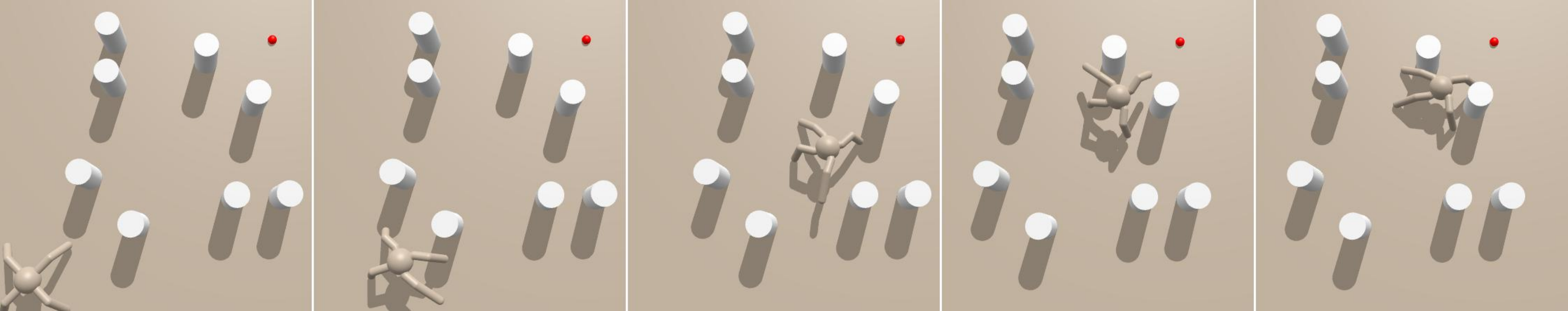}
\caption{Failure}
\label{fig:ant_local_minima}
\label{fig:ant_structure_local_minima}
\end{subfigure}
\caption{Examples of a success and a failure poor local minimum of the ant example. The obstacles are white cylinders and the target is a red sphere. Time elapses from left to right.}
\label{fig:ant_overview}
\end{figure}

\subsection{Experiments with different  $\alpha$}
This section first provides a review of dynamics and tasks, and then presents the results for each system.\\
\underline{Swimmer:} The swimmer has three links with two joints. The state $x$ consists of the position of the nose (tip of the first link), the angles of the three links, and their time derivatives. Hence, $x \in \mathbb{R}^{10}$. The control $u \in \mathbb{R}^{2}$ is torque applied to the joints. The task of the swimmer is to reach the target while avoiding spherical obstacles. We made the swimmer shorter from the original implementation so that it satisfies the constraints more easily. We let the experiments run for 300 time steps and ten times, varying the temperature parameter $\alpha$ with one environment. The experiments are considered successful when the distance between the nose of the swimmer and the target becomes less than 0.1. An example of snapshots for a successful task and that of the failed trajectory where the swimmer is captured in a poor local minimum are presented in Fig. \ref{fig:swimmer_overview}. The number of successful runs of ten at a given time step is presented in Table \ref{tab:swimmer_mpc}, where the results of the best method for the same $\alpha$ and the same time step are boldfaced. We note that the normal DDP cannot hit the target, and thus is not in the table. In most cases, SVDDP shows the best performance. \\
\underline{Ant:} The ant has four legs each of which has two joints. The state $x$ consists of the 3D position and orientation (in quaternion) of the body, angles of the legs, and its time derivatives (without quaternion). Thus, $x \in \mathbb{R}^{27}$. The control $u \in \mathbb{R}^{8}$ is torque applied to the leg joints. The task of the ant is to reach the target while keeping a prespecified distance from the center of the body to the obstacles. As in the case of the swimmer, the experiments are performed for 300 time steps and ten times. An experiment is successful when the 2D ($xy$) distance of the center of the ant's body and the target becomes less than 0.1. We present two sequences of snapshots of the two successful different trajectories in \ref{fig:ant_overview}. Fig. \ref{fig:ant_structure_local_minima} shows a poor local minimum in which the normal MPC DDP is captured and cannot escape. Table \ref{tab:ant_mpc} shows the number of successful runs out of ten as in the case of the swimmer. Here, SVDDP shows the best performance again.

\begin{table}[!h]
\footnotesize
    \caption{Number of successful runs out of ten with different $\alpha$ at 150, 200, 250, and 300 time steps.}
    \label{tab:change_alpha_swimmer_ant}
    \begin{subtable}{.49\linewidth}
\centering
\caption{Swimmer}
\label{tab:swimmer_mpc}
\resizebox{\columnwidth}{!}{
\begin{tabular}{cccccc}
\multirow{2}{*}{$\alpha$} & \multirow{2}{*}{Method} & \multicolumn{4}{c}{Time step}                                                                         \\
                       &                         & \multicolumn{1}{l}{150} & \multicolumn{1}{l}{200} & \multicolumn{1}{l}{250} & \multicolumn{1}{l}{300} \\ \hline
     & SVDDP    & \textbf{3} & \textbf{5} & \textbf{7} & \textbf{9} \\
0.01 & MG-MEDDP & 1 & 1 & 1 & 2 \\
     & UG-MEDDP & 1 & 3 & 6 & \textbf{9} \\ \hline
     & SVDDP    & \textbf{3} & \textbf{5} & \textbf{8} & \textbf{9} \\
0.1  & MG-MEDDP & 1 & 3 & 7 & 8 \\
     & UG-MEDDP & 0 & 1 & 3 & 8 \\ \hline
     & SVDDP    & 2 & \textbf{4} & \textbf{6} & \textbf{8} \\
1    & MG-MEDDP & 2 & 3 & 5 & \textbf{8} \\
     & UG-MEDDP & \textbf{4} & \textbf{4} & 5 & 6 \\ \hline
     & SVDDP    & 2 & \textbf{5} & 6 & 8 \\
5    & MG-MEDDP & \textbf{3} & \textbf{5} & 5 & 8 \\
     & UG-MEDDP & 2 & \textbf{5} & \textbf{9} & \textbf{9} \\ \hline
     & SVDDP    & \textbf{4} & \textbf{5} & 5 & \textbf{8} \\
10   & MG-MEDDP & 0 & 1 & 2 & 6 \\
     & UG-MEDDP & 1 & 4 & \textbf{7} & \textbf{8} \\ \hline
     & SVDDP    & \textbf{5} & \textbf{7} & \textbf{7} & \textbf{8} \\
20   & MG-MEDDP & 1 & 4 & 6 & 7 \\
     & UG-MEDDP & 3 & 6 & 6 & 6
\end{tabular}}
\end{subtable}%
\begin{subtable}{.49\linewidth}
\centering
\caption{Ant}
\label{tab:ant_mpc}
\resizebox{\columnwidth}{!}{
\begin{tabular}{cccccc}
\multirow{2}{*}{$\alpha$} & \multirow{2}{*}{Method} & \multicolumn{4}{c}{Time step}                                                                         \\
                       &                         & \multicolumn{1}{l}{150} & \multicolumn{1}{l}{200} & \multicolumn{1}{l}{250}& \multicolumn{1}{l}{300} \\ \hline
     & SVDDP    & \textbf{2} & \textbf{6} & \textbf{8} & 9  \\
0.01 & MG-MEDDP & 1 & 4 & \textbf{8} & \textbf{10} \\
     & UG-MEDDP & 0 & 1 & 2 & 5  \\ \hline
     & SVDDP    & 2 & 3 & \textbf{8} & \textbf{9}  \\
0.1  & MG-MEDDP & 1 & 3 & 7 & 7  \\
     & UG-MEDDP & \textbf{3} & \textbf{4} & 4 & 7  \\ \hline
     & SVDDP    & \textbf{3} & \textbf{5} & \textbf{7} & \textbf{8}  \\
1    & MG-MEDDP & 1 & 3 & \textbf{7} & 7  \\
     & UG-MEDDP & 2 & 3 & 5 & 5  \\ \hline
     & SVDDP    & \textbf{3} & 4 & \textbf{7} & 7  \\
5    & MG-MEDDP & \textbf{3} & \textbf{5} & 6 & \textbf{9}  \\
     & UG-MEDDP & 0 & 3 & 4 & 5  \\ \hline
     & SVDDP    & \textbf{3} & \textbf{5} & \textbf{6} & 7  \\
10   & MG-MEDDP & 1 & 3 & 5 & 5  \\
     & UG-MEDDP & 0 & 3 & 5 & \textbf{8}  \\ \hline
     & SVDDP    & 1 & \textbf{7} & \textbf{8} & \textbf{8}  \\
20   & MG-MEDDP & \textbf{2} & 5 & 7 & 7  \\
     & UG-MEDDP & 0 & 2 & 4 & 6 
\end{tabular}}
\end{subtable} 
\end{table}

\begin{figure}
\begin{subfigure}[b]{1\linewidth}
\includegraphics[width=\textwidth]{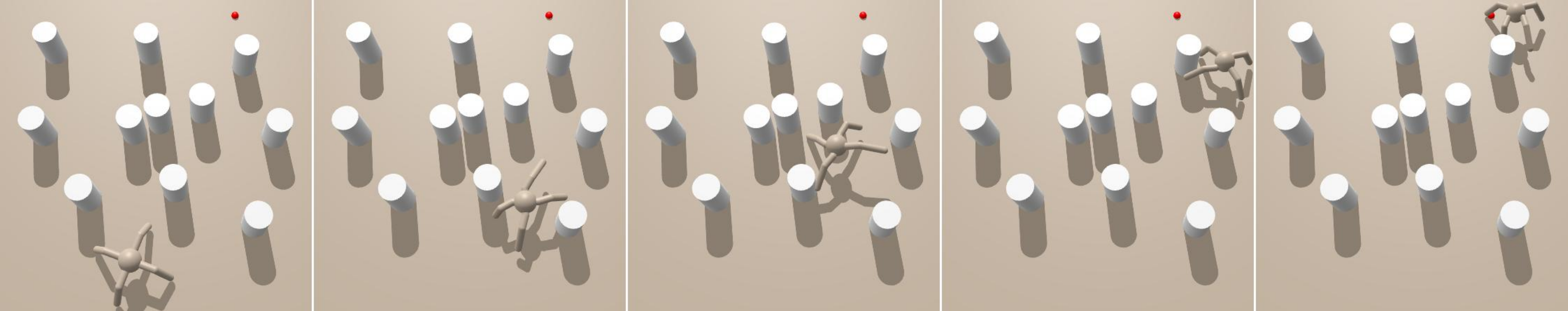}
\caption{Success}
\label{fig:ant_over_3}
\end{subfigure}
\begin{subfigure}[b]{1\linewidth}
\centering
\includegraphics[width=\textwidth]{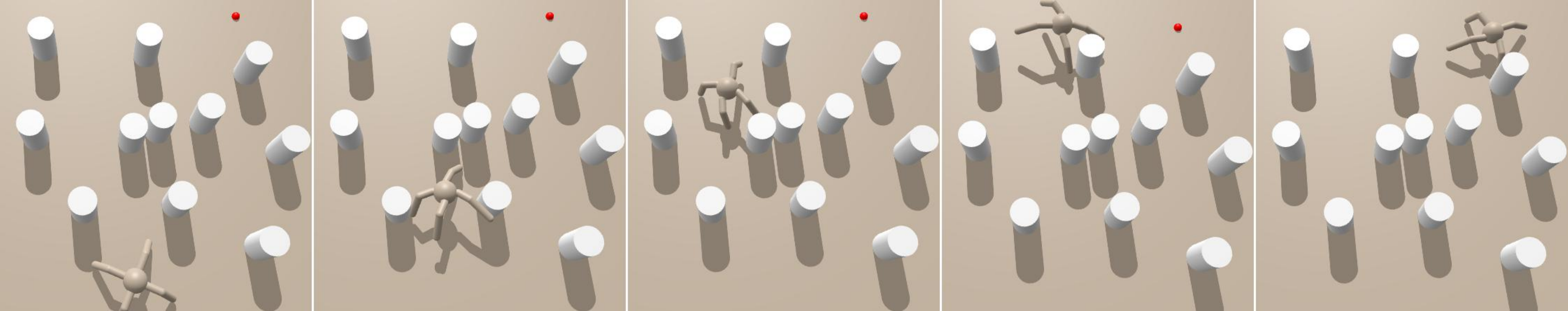}
\caption{Success}
\label{fig:ant_over_4}
\end{subfigure}
\begin{subfigure}[b]{1\linewidth}
\centering
\includegraphics[width=\textwidth]{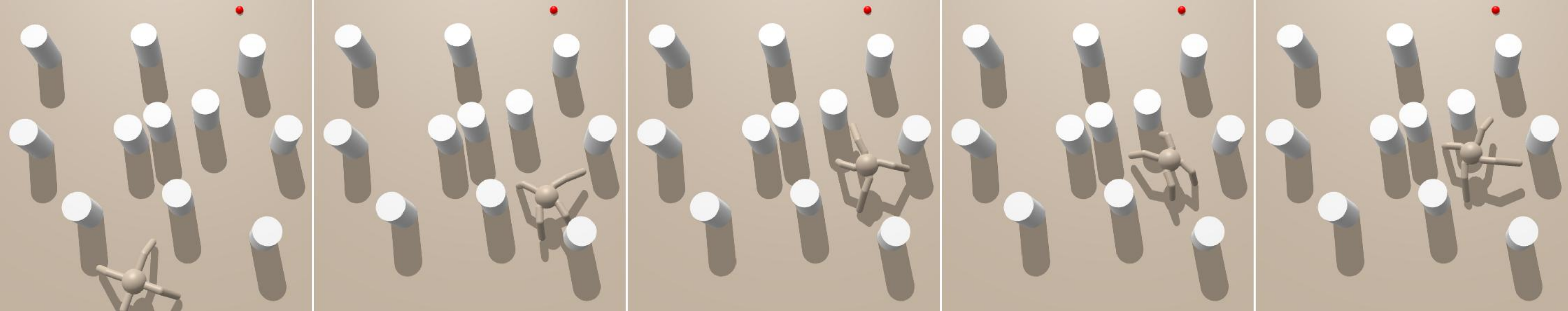}
\caption{Failure}
\label{fig:ant_over_5}
\end{subfigure}
\caption{Examples of a success and a failure poor local minima with one of the obstacle fields.}
\label{fig:ant_field_example}
\end{figure}

\subsection{Experiments of ants in obstacle fields.}
We use the ant and solve tasks similar to those in the previous section in four additional different obstacle fields. In this experiment, we fix $\alpha = 0.1$. The experiment is performed ten times in five environments, i.e., one in the previous section plus four new ones, giving 50 runs. One of the environments and two successful and one failed runs are presented in Fig. \ref{fig:ant_field_example}. The percentage of successful runs up to 150, 200, 250, 300 time steps are given in Table \ref{tab:ant_fileds}. Although MG-MEDDP performs better in the early timesteps, SVDDP eventually achieves a better success rate, showing its ability to escape poor local minima.

\begin{table}
\caption{Results of ant in five different obstacle fields. Percentage of success out of 50 runs (five environments with ten trials) for SV and MEDDPs and five runs in DDP. The best value for a time step is boldfaced.}
\label{tab:ant_fileds}
\centering
\begin{tabular}{ccccc}
\multirow{2}{*}{Method} & \multicolumn{4}{c}{Time step} \\
                        & 150   & 200   & 250   & 300   \\ \hline
SVDDP                   & 32    & \textbf{54}    & \textbf{68}    & \textbf{76 }   \\
MG-MEDDP                & \textbf{38}    & \textbf{54}    & 62    & 62    \\
UG-MEDDP                & 18    & 36    & 48    & 60    \\
DDP                     & 20    & 40    & 40    & 40   
\end{tabular}
\end{table}

\end{document}